\documentclass[11pt,a4paper]{amsart}
\usepackage[T1]{fontenc}
\usepackage[utf8]{inputenc}
\usepackage{lmodern,microtype}
\usepackage[a4paper,margin=29mm]{geometry}
\usepackage{amsmath,amssymb,amsthm,mathtools}
\usepackage{enumitem,needspace}
\usepackage[hidelinks]{hyperref}
\hypersetup{pdftitle={Nef complexity of klt Fano varieties: the full spectrum and classification below two},pdfauthor={Kiwamu Watanabe}}
\numberwithin{equation}{section}
\newtheoremstyle{plainroman}{3pt}{3pt}{\normalfont}{}{\bfseries}{.}{.5em}{}
\theoremstyle{plainroman}
\newtheorem{theorem}{Theorem}[section]
\newtheorem{proposition}[theorem]{Proposition}
\newtheorem{lemma}[theorem]{Lemma}
\newtheorem{corollary}[theorem]{Corollary}
\theoremstyle{definition}
\newtheorem{definition}[theorem]{Definition}
\theoremstyle{remark}
\newtheorem{remark}[theorem]{Remark}
\theoremstyle{plainroman}
\newtheorem*{theorem*}{Theorem}
\newtheorem*{lemma*}{Lemma}
\theoremstyle{definition}
\newtheorem*{definition*}{Definition}
\theoremstyle{remark}
\newtheorem*{remark*}{Remark}
\newcommand{\PP}{\mathbb P}
\newcommand{\QQ}{\mathbb Q}
\newcommand{\RR}{\mathbb R}
\newcommand{\ZZ}{\mathbb Z}
\newcommand{\CC}{\mathbb C}
\newcommand{\OO}{\mathcal O}
\newcommand{\Nef}{\operatorname{Nef}}
\newcommand{\Pic}{\operatorname{Pic}}
\newcommand{\Bl}{\operatorname{Bl}}
\newcommand{\Sing}{\operatorname{Sing}}
\newcommand{\eps}{\varepsilon}
\newcommand{\calS}{\mathcal S}
\newcommand{\ee}{e}

\newcommand{\vol}{\operatorname{vol}}
\newcommand{\Supp}{\operatorname{Supp}}
\newcommand{\ord}{\operatorname{ord}}
\setlist[enumerate]{itemsep=3pt,topsep=5pt}
\title[Nef complexity of klt Fano varieties]{Nef complexity of klt Fano varieties:\\ the full spectrum and classification below two}
\author{Kiwamu Watanabe}
\address{Department of Mathematics, School of Science and Engineering, Chuo University\\
1-13-27 Kasuga, Bunkyo-ku, Tokyo 112-8551, Japan}
\email{watanabe@math.chuo-u.ac.jp}
\date{}
\begin{document}
\begin{abstract}
Let $X$ be a complex klt Fano variety.  We study the rational total index
$\tau_X$ and the nef complexity
$c_X=\dim X+\rho(X)-\tau_X$.  Previous joint work of the author with
Enwright, Filipazzi, Gongyo, Moraga, Svaldi, and Wang established the
corresponding nonnegativity and zero-case characterization in the broader
framework of a special class of generalized log Calabi--Yau pairs.  Here we
give a direct proof in the Fano setting, prove that there is no value
strictly between zero and one, and classify all klt Fano varieties with
$c_X<2$, without assuming $\QQ$-factoriality.  At complexity one, after
splitting off projective-space factors, the remaining factor is a normal
quadric of Picard number one, a linear blow-up of projective space, or an
integral $(1,1)$ hypersurface in a product of projective spaces.  If
$1<c_X<2$, only two weighted-projective families occur; in particular, no
smooth Fano variety has nef complexity strictly between one and two.  We
also determine the spectrum when the dimension is allowed to vary:
\[
 \{0\}\cup\left\{2-\frac2q \ \middle| \ q\in\ZZ_{\geq2}\right\}
 \cup(\QQ\cap[2,\infty)).
\]
Most proofs in this paper were developed from arguments generated by
OpenAI's ChatGPT through discussions initiated and guided by the author;
see the \hyperref[sec:ai-disclosure]{statement on the use of generative AI}.
\end{abstract}
\maketitle
\vspace{-1.10em}

\section{Introduction}\label{sec:introduction}

All varieties are defined over $\CC$, and Fano varieties and underlying
varieties of log Fano pairs are assumed positive-dimensional.  A klt Fano
variety is a normal projective variety $X$ such that $(X,0)$ is klt and
$-K_X$ is ample and $\QQ$-Cartier; we do not assume $\QQ$-factoriality.

The Kobayashi--Ochiai theorem \cite{KobayashiOchiai} bounds the Fano index
of a smooth Fano $n$-fold by $n+1$, with equality characterizing projective
space.  Mukai \cite{Mukai} proposed a refinement involving the Picard number,
with equality characterizing products of projective spaces.  Much of the
progress on Mukai-type problems, especially on the generalized Mukai
conjecture formulated in terms of the pseudoindex, has proceeded through Mori
theory and deformation of rational curves; see, for example,
\cite{WisniewskiMukai,BCDD,NovelliOcchetta,FujitaMukai}.  These arguments use
unsplit or locally unsplit families of rational curves, estimates for their
loci, and the geometry of extremal contractions.  They are particularly
effective when the pseudoindex is large relative to the dimension or when the
extremal contractions satisfy additional restrictions.  In arbitrary
dimension, however, it is difficult to construct and control enough
numerically independent families while retaining the sharp dimension
estimates needed to recover the Picard number and the global product
structure.  This makes a direct extension of the rational-curve approach to
a uniform higher-dimensional classification difficult.

The total index suggests a complementary direction.  Gongyo \cite{Gongyo}
introduced the total index, which measures anticanonical positivity through
decompositions into nef Cartier divisors, while related complexity invariants
for log Calabi--Yau pairs were introduced by Gongyo--Moraga
\cite{GongyoMoraga}; see also \cite{BMSZ}.  Instead of following individual
rational curves and their deformations, one can study all nef directions
appearing in an anticanonical decomposition simultaneously.  The approach of
the present paper develops this idea by combining optimal nef decompositions,
vanishing and finite differences of multivariable Euler polynomials, and
mixed-intersection theory.  The support of the positive mixed intersections
carries a discrete-convex exchange structure, whose rank is identified with
the numerical dimensions of partial sums of the nef divisors and hence, in
the semiample situation, with dimensions of contraction images.  This gives
a route from anticanonical positivity to contractions and product splittings
that does not require constructing a prescribed collection of families of
rational curves.

Related developments include the generalized Mukai conjecture for
$\QQ$-factorial spherical Fano varieties \cite{GHP} and the integral and
rational Gongyo indices for smooth toric Fano varieties \cite{SatoGongyo}; in
the latter setting the rational Gongyo index is the rational total index used
below.  We use the rational version:

\begin{definition}\label{def:nef-complexity}
For a complex klt Fano variety $X$, put $n=\dim X$ and $\rho=\rho(X)$, and define
\begin{align}
 \tau_X&=\sup\left\{\sum_{i=1}^{s}a_i\ \middle|\
 \begin{array}{l}
 a_i\in\QQ_{>0},\quad L_i\text{ nef Cartier},\quad L_i\not\equiv0,\\[-2pt]
 -K_X\equiv\sum_{i=1}^{s}a_iL_i
 \end{array}\right\},\label{eq:definition}\\[-2pt]
 c_X&=n+\rho-\tau_X.\nonumber
\end{align}
We call $\tau_X$ the rational total index and $c_X$ the nef complexity of $X$.
\end{definition}

In the joint work \cite{EFGMSWW}, nef complexity was studied through a
special class of generalized klt log Calabi--Yau pairs.  That work established
the corresponding nonnegativity and zero-case theorem, with products of
projective spaces in the equality case; it also classified smooth Fano
threefolds of nef complexity at most one and proved a Mukai-type statement
under a fibre-type hypothesis on extremal contractions.  It did not determine
the first positive values of nef complexity in arbitrary dimension, nor the
geometry of klt Fano varieties immediately beyond the equality case.

The present paper shows that the divisor-theoretic viewpoint reaches
substantially further.  We prove that there is no value strictly between zero
and one, classify all klt Fano varieties of complexity one in arbitrary
dimension, classify the entire range $1<c_X<2$, and determine the full
spectrum when the dimension is allowed to vary.  In particular, the
method detects both fibre-type and birational phenomena, including the
linear blow-ups and the weighted-projective families appearing below, and
does so without assuming $\QQ$-factoriality.  Thus the classification below
two should also be viewed as a first test of a broader strategy: use an
optimal anticanonical decomposition and its mixed intersections as a
replacement for direct control of many rational-curve families.  The
previously known nonnegativity and zero case are included in the main theorem
for context.

The use of generative AI in developing the proofs, including the
respective contributions of ChatGPT and the author, is described
in Subsection~\ref{sec:ai-disclosure}.

For a smooth Fano variety with $\rho(X)=1$, $\tau_X$ coincides with the
Fano index, while products of projective spaces satisfy $c_X=0$.

\subsection{Main results}

The main theorem gives the classification below two in arbitrary dimension,
without assuming $\QQ$-factoriality.

\begin{theorem}\label{thm:main}
Let $X$ be a complex klt Fano variety.  Then the
following statements hold.  We write $\PP(w_0,\ldots,w_N)$ for the
weighted projective space with weights $w_0,\ldots,w_N$.
\begin{enumerate}[label=(\roman*)]
\item One has $c_X\geq0$.  Moreover, $c_X=0$ holds exactly when
\[
 X\simeq\prod_{j=1}^{t}\PP^{d_j},\qquad d_j\geq1.
\]
\item If $c_X<1$, then $c_X=0$.
\item If $c_X=1$, then
\begin{equation}\label{eq:main-product}
 X\simeq Y\times\prod_{j=1}^{t}\PP^{d_j},\qquad d_j\geq1,
\end{equation}
where the projective-space product may be empty and $Y$ is one of the
following:
\begin{enumerate}[label=(\alph*)]
\item a normal quadric $Q^d\subset\PP^{d+1}$ of dimension $d\geq2$ and
Picard number one; equivalently, $d\geq3$ and the quadratic form has rank
at least three, or $d=2$ and $Q^2$ is the quadric cone;
\item $\Bl_{\PP^\ell}\PP^d$, where $d\geq2$ and
$0\leq\ell\leq d-2$;
\item an integral hypersurface of bidegree $(1,1)$ in
$\PP^a\times\PP^b$, where $a,b\geq2$.
\end{enumerate}
Conversely, all these varieties have nef complexity one.
\item If $1<c_X<2$, then $X$ is a product as in
\eqref{eq:main-product}, where $Y$ is one of the following:
\begin{enumerate}[label=(\alph*)]
\item $\PP(1,1,\underbrace{q,\ldots,q}_{d-1})$, where $d\geq2$ and
$q\geq3$; in this case $c_X=2-2/q$;
\item $\PP(1,1,1,\underbrace{2,\ldots,2}_{d-2})$, where $d\geq3$; in
this case $c_X=3/2$.
\end{enumerate}
Conversely, the varieties in (a) and (b), and their
products with projective spaces, have the stated nef complexity.
\end{enumerate}
The hypersurfaces in (iii) (c) are normal, terminal, and
Gorenstein.
\end{theorem}

In particular, there is no value strictly between zero and one, and the
two weighted-projective families in Theorem~\ref{thm:main} determine the
positive spectrum below two.

\begin{corollary}\label{cor:spectrum}
The values of $c_X$ strictly less than two are exactly
\[
 \{0\}\cup\left\{2-\frac2q\ \middle| \ q\in\ZZ,\ q\geq2\right\}.
\]
\end{corollary}

Allowing the dimension to vary, one obtains the full spectrum.

\begin{theorem}\label{thm:full-spectrum}
Let
\[
 \operatorname{Spec}_{\mathrm{klt}}(c)
 :=\left\{c_X\ \middle|\  X\text{ is a complex klt Fano variety}\right\}.
\]
Then
\begin{equation}\label{eq:full-spectrum}
 \operatorname{Spec}_{\mathrm{klt}}(c)
 =\{0\}\cup
 \left\{2-\frac2q\ \middle| \ q\in\ZZ_{\geq2}\right\}
 \cup\bigl(\QQ\cap[2,\infty)\bigr).
\end{equation}
Moreover, every rational number at least two is realized by a
$\QQ$-factorial klt Fano variety of Picard number two.
\end{theorem}

Thus two is the first accumulation point of the spectrum, and every rational
value at least two already occurs in Picard number two.

\subsection{Proof strategy and the main mechanism}

The proof is organized around a passage from anticanonical decompositions
to finite combinatorial data and then back to the geometry of contractions.
We describe this mechanism in some detail, since it is used repeatedly in
the classification and is not specific to the threshold $c_X<2$.  No
familiarity with Lorentzian polynomials or $M$-convexity is assumed; the
precise definitions and the consequences needed below are recalled in
Section~\ref{sec:preliminaries}.

The first step is to replace the supremum defining $\tau_X$ by a finite
linear programme.  The nef Cartier monoid has a finite Hilbert basis, and a
maximizing rational vertex gives an optimal decomposition
\[
 -K_X\equiv\sum_{i=1}^s a_iL_i,
 \qquad
 a_i=m_i+1-\eps_i,
 \qquad
 0\leq\eps_i<1,
\]
whose supporting classes are linearly independent.  For
$k=(k_1,\ldots,k_s)$ with $0\leq k_i\leq m_i$, the divisor
$-K_X-\sum_i k_iL_i$ remains ample.  Kawamata--Viehweg vanishing therefore
forces the multivariable Euler polynomial
\[
 P(k)=\chi\!\left(X,\OO_X\!\left(-\sum_i k_iL_i\right)\right)
\]
to vanish at every nonzero lattice point of this rectangle, while
$P(0)=1$.  We call this finite rectangle the \emph{vanishing box}.
Applying finite-difference operators to $P$ converts this vanishing pattern
into numerical restrictions on the decomposition.  In particular,
\[
 \sum_i m_i\leq n,
 \qquad
 c_X=\left(n-\sum_i m_i\right)+(\rho-s)+\sum_i\eps_i.
\]
Thus the nef complexity is separated into two integral defects and one
fractional defect.  When $c_X$ is small, all three defects are forced to be
small simultaneously; this is the first source of rigidity.

The second step is to organize the positive mixed intersections of the
$L_i$.  Consider the volume polynomial
\[
 \left(\sum_i t_iL_i\right)^n
\]
and the set
\[
 \calS=\left\{\alpha\in\ZZ_{\geq0}^s\ \middle|\
 |\alpha|=n,\ L^\alpha>0\right\}
\]
of exponent vectors of its positive coefficients.  Br\"and\'en--Huh's
Lorentzian theory \cite{BrandenHuh} implies that $\calS$ is $M$-convex.
For the present argument, the important content is the symmetric exchange
property: positivity of two mixed intersections can be transported by
moving one unit of exponent from one nef direction to another.  Moreover,
for every $J\subseteq\{1,\ldots,s\}$, the associated rank function is
\[
 \max_{\alpha\in\calS}\sum_{i\in J}\alpha_i
 =\nu\!\left(\sum_{i\in J}L_i\right).
\]
Hence a discrete-convex rank attached to mixed intersections is identified
with the numerical dimension of a partial sum of nef divisors; when that
sum is semiample, this is the dimension of the image of its contraction.
This gives the basic bridge from mixed-intersection combinatorics back to
birational geometry.

The finite-difference identities and the exchange property are then used
together.  Finite differences produce specific vanishing and nonvanishing
mixed intersections, while $M$-convexity propagates positivity among nearby
exponent vectors.  The resulting rank constraints separate the nef
directions into \emph{active} directions, which carry the residual mixed
intersection, and \emph{inactive} directions.  The latter have sufficiently
rigid contraction behaviour to split off as projective-space factors.
After these factors are removed, the classification problem is reduced to
a factor of small Picard number carrying only the active directions.

At complexity one this leaves a factor of Picard number at most two, leading
to the quadric, linear blow-up, and $(1,1)$-hypersurface cases.  In the range
$1<c_X<2$, the same mechanism yields a weighted mixed-intersection identity
that permits at most two active directions.  A birational degree lemma and
a surface argument exclude the two-active case, while the one-active case
reduces to a Picard-number-one polarized variety and gives the two
weighted-projective families in Theorem~\ref{thm:main}.  The full spectrum
is then completed by an explicit Picard-number-two construction realizing
every rational value at least two.

The individual ingredients above come from convex geometry, vanishing
theorems, intersection theory, and discrete convexity.  The point of the
argument is their combination into a single procedure
\[
\begin{aligned}
 \text{optimal nef decomposition}
 &\longrightarrow \text{Euler-polynomial vanishing and finite differences}\\
 &\longrightarrow \text{mixed-intersection exchange and rank}\\
 &\longrightarrow \text{contractions and product splittings}.
\end{aligned}
\]
This procedure is useful beyond the final list of varieties obtained here:
it provides a way to turn a sufficiently positive anticanonical
decomposition into combinatorial constraints on mixed intersections and
then recover geometric information from those constraints.

\subsection{Organization}

Sections \ref{sec:preliminaries}--\ref{sec:equality} develop the basic
nef-cone, finite-difference, and mixed-intersection tools, including the gap
below one.  Section~\ref{frac:sec:full} treats the remaining full-degree
case, and Section~\ref{frac:sec:mixed} establishes the inactive-factor
reduction for $1\leq c_X<2$.  Section~\ref{sec:one} treats complexity one.
Sections \ref{frac:sec:birational}--\ref{frac:sec:classification} handle
the range $1<c_X<2$.  Section~\ref{sec:converse} verifies the classified
varieties, and Section~\ref{sec:spectrum} determines the full spectrum.

\medskip
\noindent\textbf{Note added.}
During the final checks of this manuscript, OpenAI released the preprint
\cite{OpenAIMukai}, which claims a proof of the generalized Mukai conjecture
for smooth complex Fano varieties.

\subsection{Statement on the use of generative AI}\label{sec:ai-disclosure}
The research leading to this paper has its origins in Gongyo's work on
the total index and the Mukai-type conjecture \cite{Gongyo}, as well as
related work by Gongyo and Moraga on generalized complexity
\cite{GongyoMoraga}. At the AIM workshop
\emph{Higher-dimensional log Calabi--Yau pairs}, held in
September--October 2024, Gongyo suggested that he investigate the
Mukai-type conjecture while the author study the (generalized) Mukai
conjecture. This initiated the collaboration that eventually led to
\cite{EFGMSWW}.

In that joint work, we proved the non-negativity of nef complexity and
characterized its zero case by products of projective spaces, using
results for generalized log Calabi--Yau pairs. We additionally sought
to establish that $c_X<1$ implies $c_X=0$ for smooth Fano varieties and
to classify those with $c_X=1$. However, we were able to complete these
classifications only in dimension three.

The present research therefore began with the goal of extending these
threefold results to smooth Fano varieties of arbitrary dimension,
initially using the theory of generalized log Calabi--Yau pairs. During
the spring and summer of 2026, particularly from approximately April to
July, the author investigated this approach, but attempts to carry out
the necessary inductive arguments did not succeed.

In parallel, the author had been studying the total index and nef
complexity of smooth Fano threefolds. In early 2026, the author asked
ChatGPT to compute these invariants for the known families. During
these discussions, ChatGPT proposed using the Hilbert basis of the nef
Cartier monoid and linear programming to compute the total index. This
computational approach was not pursued further at that time as a method
for proving general classification results.

When the author returned to the threefold computations in early
September 2026, the Hilbert-basis and linear-programming approach again
became central. The author then asked ChatGPT to apply this method to
obtain an alternative proof of the non-negativity of nef complexity and
the characterization of its zero case for klt Fano varieties, without
relying on the theory of generalized log Calabi--Yau pairs.

After several rounds of discussion, ChatGPT developed an alternative
argument. The author recognized that this approach might also
establish the stronger statement that $c_X<1$ implies $c_X=0$. Further
discussions led to the classification of the case $c_X=1$, followed by
the classification of varieties satisfying $1<c_X<2$. These results
were obtained through repeated interactions in which arguments were
proposed, examined, corrected, and refined.

After obtaining the classification below two, the author became
interested in the possible values of nef complexity for $c_X\ge2$.
When the author raised this question, ChatGPT proposed, in its first
response, essentially the construction used in Proposition~11.1,
together with the formula
\[
c_X=2+\frac{2(d-1)}{q},
\]
which realizes every rational number at least two. The details of the
construction and its proof were subsequently examined and refined.

Most of the proofs in this paper originated from arguments generated by
ChatGPT, and many of the central technical ideas were also proposed by
the model. In particular, ChatGPT played a major role in developing the
arguments involving Hilbert-basis linear programming, vanishing and
finite differences of Euler-characteristic polynomials,
mixed-intersection theory, Lorentzian polynomials, and $M$-convexity, as
well as the classification arguments.

The author's contributions included formulating the research problems,
recognizing connections with previous work, identifying potential
applications of the methods proposed by ChatGPT, directing the
subsequent investigations, examining mathematical arguments, testing
examples and possible counterexamples, and identifying gaps or errors
requiring correction. The transition from the threefold computations
to the general study of nef complexity was prompted by the author's
decision to apply the computational approach to the problems left open
in the earlier joint work.

Most of the principal results were obtained in early September 2026,
primarily using OpenAI's GPT-5.6 Sol. GPT-6 Sol was used subsequently,
mainly for examining and revising proofs and preparing the manuscript.
Earlier investigations had also involved ChatGPT, using models
available at those times.

The verification and revision process was itself partly AI-assisted.
The author does not claim to have independently derived every argument
without AI assistance. Nevertheless, the author has personally checked
all proofs in the final manuscript. The author determined the final
mathematical statements and presentation and assumes full responsibility
for the mathematical content, including its correctness, originality,
and attribution.

\section{Preliminaries on nef Cartier divisors}\label{sec:preliminaries}

We use $\equiv$ for numerical equivalence and $\sim$ for linear equivalence.
A \emph{contraction} is a projective surjective morphism to a normal variety
with connected fibres.  All boundaries considered below are effective
$\QQ$-divisors.  We recall only the convex-geometric, discrete-convex, and
intersection-theoretic facts that will be used later.

\subsection{Lattices, rational polyhedra, polyhedral cones, and Hilbert bases}

We use the standard terminology for rational polyhedra, lattices, faces,
vertices, and active inequalities from Bruns--Gubeladze
\cite[Chapters~1--2]{BG}.  Thus a lattice $M\subset V$ induces the rational
structure $V_{\QQ}=M\otimes_{\ZZ}\QQ$, and ``rational'' always refers to the
chosen rational structures.  We write $\operatorname{Conv}(S)$ for convex hulls.
Whenever $\RR^r$ or $\ZZ^r$ is used, $\ee_i$ denotes the $i$-th standard unit
vector.

\begin{lemma}\label{lem:polyhedral-basic}
Let $P$ be a polyhedron in a finite-dimensional real vector space equipped
with a rational structure.
\begin{enumerate}[label=(\roman*)]
\item If $P$ is rational, every vertex of $P$ is rational.
\item If $P$ is nonempty and bounded, then it is the convex hull of its finitely
many vertices.
\item If a linear functional $\ell$ attains its maximum $T$ on $P$, then
$\{x\in P\mid\ell(x)=T\}$ is a face.  In particular, on a nonempty bounded
rational polyhedron, $\ell$ attains its maximum at a rational vertex.
\item Let $B:\RR^N\to V$ be linear and
$P_b:=\{x\in\RR_{\geq0}^N\mid B(x)=b\}$.  If $B$ is rational and
$b\in V_{\QQ}$, then $P_b$ is rational.  If $x$ is a vertex and
$I=\{i\mid x_i>0\}$, the vectors $B(\ee_i)$, $i\in I$, are linearly
independent.
\end{enumerate}
\end{lemma}
\begin{proof}
Assertions (i)--(iii) are standard facts on rational
polyhedra; see \cite[Chapter~1]{BG}, in particular
\cite[Proposition~1.69~(c), Theorems~1.10~(b) and~1.26]{BG}.
For (iv), rationality is immediate. If the vectors $B(\ee_i)$ with
$x_i>0$ were linearly dependent, there would be a nonzero vector $u$
supported on these indices with $B(u)=0$. For sufficiently small
$\varepsilon>0$, the distinct points $x\pm\varepsilon u$ would both lie in
$P_b$ and have midpoint $x$, contradicting that $x$ is a vertex.
\end{proof}

If $V$ has lattice $M$, a finite set $S\subset M$ is \emph{hole-free} if
$S=\operatorname{Conv}(S)\cap M$.  We also use the standard facts that a
rational polyhedral cone has rational polyhedral dual and, by Gordan's lemma,
$C\cap M$ is finitely generated; for a pointed cone its unique minimal
additive generating set is the \emph{Hilbert basis}.  See
\cite[Theorem~1.15, Proposition~1.69, Lemma~2.9, Definition~2.15 and
Proposition~2.16]{BG}.

\subsection{\texorpdfstring{$M$}{M}-convex sets and discrete polymatroids}
\label{subsec:Mconvex-background}

For an integer $s\geq1$, let $[s]=\{1,\ldots,s\}$ and $2^{[s]}$ the set of all subsets of $[s]$.  For
$x=(x_1,\ldots,x_s)\in\RR^s$ and $T\subseteq[s]$, write
\(x(T):=\sum_{i\in T}x_i\) and \(|x|:=x([s])\).
We use the componentwise order on $\RR^s$.  
A nonempty finite set $P\subseteq\ZZ_{\geq0}^s$ is a \emph{discrete
polymatroid} if it is downward closed and satisfies the augmentation axiom:
if $u,v\in P$ and $|u|<|v|$, then $u+\ee_i\in P$ for some index $i$ with
$u_i<v_i$.  A \emph{base} is a maximal element for the componentwise order.
These are the conventions of \cite{HH}.

A nonempty set $S\subseteq\left\{x\in\ZZ_{\geq0}^s\ \middle|\ |x|=d\right\}$ is \emph{$M$-convex} if, for every
$\alpha,\beta\in S$ and every $i$ with $\alpha_i>\beta_i$, there is an index
$j$ with $\alpha_j<\beta_j$ such that
\begin{equation}\label{eq:M-convex-exchange}
 \alpha-\ee_i+\ee_j\in S,
 \qquad
 \beta+\ee_i-\ee_j\in S.
\end{equation}
This is the symmetric exchange axiom; see \cite[Section~4.1]{Murota}.

A set function $r:2^{[s]}\to\RR$ is \emph{submodular} if
$r(T)+r(U)\geq r(T\cup U)+r(T\cap U)$ for all $T,U\subseteq[s]$.
For an $M$-convex set $S\subseteq\{x\in\ZZ_{\geq0}^s\mid |x|=d\}$, its
\emph{rank function} is $r_S(T):=\max_{\alpha\in S}\alpha(T)$.

\begin{proposition}\label{prop:Mconvex-standard}
Let $S\subseteq\left\{x\in\ZZ_{\geq0}^s\ \middle|\ |x|=d\right\}$ be a nonempty finite $M$-convex set.  Then:
\begin{enumerate}[label=(\roman*)]
\item $r_S$ is submodular.
\item For $u\in\ZZ_{\geq0}^s$, there exists $\alpha\in S$ with
$u\leq\alpha$ if and only if
\begin{equation}\label{eq:Mconvex-extension}
 u(T)\leq r_S(T)\quad\text{for every }T\subseteq[s].
\end{equation}
\item The set $S$ is hole-free:
\begin{equation}\label{eq:Mconvex-hole-free}
 S=\operatorname{Conv}(S)\cap\ZZ^s.
\end{equation}
\end{enumerate}
In particular, if $|u|=d$, then $u\in S$ if and only if
$u(T)\leq r_S(T)$ for every $T\subseteq[s]$.
\end{proposition}
\begin{proof}
Let $P_S$ be the downward closure of $S$.  By
\cite[Theorem~2.3, Proposition~1.2 and Theorem~3.4]{HH}, $P_S$ is the
discrete polymatroid with base set $S$ and rank function $r_S$, and its
integral points are exactly those satisfying $x(T)\leq r_S(T)$ for all
$T\subseteq[s]$.  This gives (i) and (ii); restricting to
$|x|=d$ gives (iii) and the final assertion.
\end{proof}

\subsection{The nef cone and descent of Cartier divisors}

For a normal projective variety $X$, put
\(\Lambda_X:=\Pic(X)/\!\equiv_{\mathrm{num}},\)
viewed as the lattice of Cartier numerical classes in $N^1(X)_{\RR}$.
The nef cone is pointed in $N^1(X)_{\RR}$.

\begin{lemma}\label{lem:nef-monoid}
Let $(X,\Delta)$ be a projective klt log Fano pair.  Then $\Nef(X)$ is
rational polyhedral, every nef Cartier divisor on $X$ is semiample, and the
monoid $\Nef(X)\cap\Lambda_X$ is finitely generated.  Hence it has a finite
Hilbert basis.
\end{lemma}
\begin{proof}
The Cone Theorem gives rational polyhedrality of the Mori cone, hence of its
dual $\Nef(X)$; see \cite[Theorem~3.7]{KM}.  If $D$ is nef Cartier, then
$D-(K_X+\Delta)$ is ample, so the base-point-free theorem
\cite[Theorem~3.24]{KM} makes $D$ semiample.  Finite generation and the
existence of the Hilbert basis follow from Gordan's lemma.
\end{proof}

The following lemma packages the contraction theorem and relative vanishing
in the form used below.

\begin{lemma}\label{lem:descent}
Let $(X,\Delta)$ be a projective klt pair and let $f:X\to Y$ be a contraction
such that $-(K_X+\Delta)$ is $f$-ample.  If a Cartier divisor $D$ on $X$ is
numerically trivial over $Y$, then
\(D\sim f^*H\)
for a Cartier divisor $H$ on $Y$.  If $f$ is the contraction defined by the
semiample divisor $D$, then $H$ is ample.  Moreover,
\begin{equation}\label{eq:relative-vanishing}
 R^qf_*\OO_X=0\qquad(q>0),
\end{equation}
and, for every Cartier divisor $B$ on $Y$,
\begin{equation}\label{eq:euler-descent}
 \chi\bigl(X,f^*\OO_Y(B)\bigr)=\chi\bigl(Y,\OO_Y(B)\bigr).
\end{equation}
\end{lemma}
\begin{proof}
The descent statement and $f_*\OO_X=\OO_Y$ are exactly the klt case of
Fujino's contraction theorem \cite[Theorem~16.4~(2),~(3)]{FujinoFundamental};
see also \cite[Lemma~2.9]{EFGMSWW}. If $f$ is defined by $D$, the descended
divisor is ample. Relative Kawamata--Viehweg vanishing gives
\eqref{eq:relative-vanishing} by \cite[Corollary~5.7.7]{FujinoMMP}, and the
projection formula and Leray give \eqref{eq:euler-descent}.
\end{proof}

In particular, on a projective klt log Fano pair, numerical equivalence of
$\QQ$-Cartier divisors implies $\QQ$-linear equivalence.

\begin{lemma}\label{lem:ample-sum}
Let $(X,\Delta)$ be a projective klt log Fano pair of
dimension $n$.  Suppose
\[
 -(K_X+\Delta)\equiv\sum_{i=1}^s a_iL_i,
 \qquad a_i>0,
\]
where the $L_i$ are nef $\RR$-Cartier divisors.  Then every positive linear
combination $\sum_i b_iL_i$, $b_i>0$, is ample.  In particular,
$A:=\sum_iL_i$ is ample.  Moreover, if $D$ is a nonzero nef numerical
Cartier class, then $D\cdot A^{n-1}>0$.
\end{lemma}
\begin{proof}
Choose $\delta>0$ with $b_i-\delta a_i\geq0$ for all $i$.  Then
$\sum_i b_iL_i\equiv\delta(-(K_X+\Delta))+\sum_i(b_i-\delta a_i)L_i$ is
ample.  The final assertion follows from semiampleness of a nef Cartier
representative.
\end{proof}

\begin{lemma}\label{lem:product-index}
Let $X$ be a complex klt Fano variety and let $r\geq1$. Then
\[
 \tau_{X\times\PP^r}=\tau_X+r+1,
 \qquad
 c_{X\times\PP^r}=c_X.
\]
\end{lemma}
\begin{proof}
Let $p:X\times\PP^r\to X$ and $q:X\times\PP^r\to\PP^r$ be the
projections, and put $H=q^*\OO_{\PP^r}(1)$.  We have
\[
 \Pic(X\times\PP^r)=p^*\Pic(X)\oplus\ZZ H,
\]
and a Cartier divisor $p^*D+bH$ is nef if and only if $D$ is nef and
$b\geq0$.

Consider a nef Cartier decomposition
\[
 -K_{X\times\PP^r}\equiv\sum_i a_i(p^*D_i+b_iH),
 \qquad a_i>0.
\]
For the indices with $D_i\not\equiv0$, restriction to
$X\times\{\mathrm{pt}\}$ gives a nef Cartier decomposition of $-K_X$,
so the sum of the corresponding coefficients is at most $\tau_X$.
For the remaining indices, $D_i\equiv0$ and hence $b_i\geq1$; since the
total $H$-coefficient is $r+1$, their coefficient sum is at most
$r+1$.  Thus
\(\tau_{X\times\PP^r}\leq\tau_X+r+1\).
Conversely, adjoining $(r+1)H$ to any nef Cartier decomposition of
$-K_X$ increases its coefficient sum by $r+1$.  Taking the supremum over
all such decompositions gives the reverse inequality. Therefore
\[
 \tau_{X\times\PP^r}=\tau_X+r+1.
\]
Since $\dim(X\times\PP^r)=\dim X+r$ and
$\rho(X\times\PP^r)=\rho(X)+1$, the equality
$c_{X\times\PP^r}=c_X$ follows.
\end{proof}

\subsection{Numerical dimension}

For a nef $\RR$-Cartier divisor $D$ on an $n$-dimensional projective variety
$X$, we use the numerical dimension
\begin{equation}\label{eq:numerical-dimension-definition}
 \nu(D):=\max\left\{k\in\{0,\ldots,n\}\ \middle|\  D^kA^{n-k}>0\right\},
\end{equation}
where $A$ is ample.  This is the usual numerical dimension of a nef class;
compare \cite[Remark~2.3.17]{Lazarsfeld}.

\begin{lemma}\label{lem:numdim-basic}
The number in \eqref{eq:numerical-dimension-definition} is independent of
$A$. If $D$ is semiample and $f:X\to Y$ is its contraction, then
$\nu(D)=\dim Y$; see \cite[Remark~2.3.17]{Lazarsfeld}.
\end{lemma}

\subsection{Mixed intersections and the volume polynomial}

Let $L_1,\ldots,L_s$ be nonzero nef Cartier divisors on an
$n$-dimensional projective variety $X$ such that
$A=L_1+\cdots+L_s$ is ample.  For $\alpha\in\ZZ_{\geq0}^s$ write
$\alpha!=\prod_i\alpha_i!$ and $L^\alpha=\prod_iL_i^{\alpha_i}$.  Define
\begin{equation}\label{eq:volume-polynomial}
 \operatorname{vol}_L(t)
 =\left(\sum_{i=1}^s t_iL_i\right)^n
 =\sum_{|\alpha|=n}\frac{n!}{\alpha!}L^\alpha t^\alpha
\end{equation}
and
\begin{equation}\label{eq:support}
 \calS:=\left\{\alpha\in\ZZ_{\geq0}^s\ \middle|\ |\alpha|=n,\ L^\alpha>0\right\}.
\end{equation}
By Br\"and\'en--Huh \cite[Theorem~4.6 and Corollary~4.7]{BrandenHuh}, the
support $\calS$ is $M$-convex.

We also use the mixed Khovanskii--Teissier inequality
\begin{equation}\label{eq:mixed-KT}
 (BCN_1\cdots N_{n-2})^2
 \geq
 (B^2N_1\cdots N_{n-2})(C^2N_1\cdots N_{n-2}),
\end{equation}
for nef classes $B,C,N_1,\ldots,N_{n-2}$; this is the case $p=2$ of
\cite[Variant~1.6.2]{Lazarsfeld}.

\begin{lemma}\label{lem:exchange}
Let $L_1,\ldots,L_s$ and $\calS$ be as above.  Then $\calS$ has the
symmetric exchange property \eqref{eq:M-convex-exchange}.  For every
$J\subseteq[s]$,
\begin{equation}\label{eq:numdim}
 \nu\left(\sum_{i\in J}L_i\right)
 =\max_{\alpha\in\calS}\sum_{i\in J}\alpha_i.
\end{equation}
Moreover, fix distinct $p,q\in[s]$ and nonnegative integers $\gamma_i$ for
$i\neq p,q$ with $\sum_{i\neq p,q}\gamma_i\leq n$.  Put
\(r:=n-\sum_{i\neq p,q}\gamma_i\)
and, for $0\leq k\leq r$, set
\(b_k:=L_p^kL_q^{r-k}\prod_{i\neq p,q}L_i^{\gamma_i}.\)
Then $b_k^2\geq b_{k+1}b_{k-1}$ for $1\leq k\leq r-1$, and
$\{k\mid b_k>0\}$ is empty or an interval of consecutive integers.
\end{lemma}
\begin{proof}
The symmetric exchange property is the $M$-convexity of $\calS$.  Put
$D_J=\sum_{i\in J}L_i$.  For $0\leq k\leq n$, multinomial expansion gives
\begin{equation}\label{eq:numdim-expansion}
 D_J^kA^{n-k}
 =\sum_{\substack{|\beta|=k,\ \beta_i=0\ (i\notin J)\\
                   |\gamma|=n-k}}
   \frac{k!}{\beta!}\frac{(n-k)!}{\gamma!}
   L^{\beta+\gamma}.
\end{equation}
All terms are nonnegative.  The right-hand side is positive exactly when
there is $\alpha\in\calS$ with $\alpha(J)\geq k$, which proves
\eqref{eq:numdim}.

For log-concavity, apply \eqref{eq:mixed-KT} with $B=L_p$, $C=L_q$, and
\[
 N_1\cdots N_{n-2}=L_p^{k-1}L_q^{r-k-1}\prod_{i\neq p,q}L_i^{\gamma_i}.
\]
The two factors on the right are $b_{k+1}$ and $b_{k-1}$.

Finally, if $b_a,b_c>0$ with $a<c$, then the corresponding exponent vectors
$\alpha^{(a)},\alpha^{(c)}\in\calS$ have every integral vector
$\alpha^{(k)}$, $a\leq k\leq c$, on the line segment joining them.  Since
$\calS$ is hole-free by Proposition~\ref{prop:Mconvex-standard}, one has
$\alpha^{(k)}\in\calS$ and hence $b_k>0$ for all $a\leq k\leq c$.
\end{proof}

\begin{remark}\label{rem:geom-polymatroid}
If the $L_i$ are semiample, choose $r>0$ so that every $rL_i$ is base-point
free and let $\varphi:X\to Z\subset\prod_iY_i$ be the product of the
associated contractions.  Since $\sum_iL_i$ is ample, $\varphi$ is finite;
moreover, the image of each partial projection has dimension
$\nu(\sum_{i\in T}L_i)$, and its positive multidegrees are exactly the
$L^\alpha>0$, up to a common positive scalar.  Thus \eqref{eq:numdim} is the
usual rank formula for the discrete polymatroid of multidegrees; compare
\cite[Theorem~A]{CastilloEtAl}.
\end{remark}

For $T\subseteq[s]$, set
\begin{equation}\label{eq:rank}
 r(T):=\max_{\alpha\in\calS}\alpha(T)
      =\nu\left(\sum_{i\in T}L_i\right).
\end{equation}
By Proposition~\ref{prop:Mconvex-standard}, $r$ is submodular and is the
rank function of the discrete polymatroid determined by $\calS$.  Moreover,
for $u\in\ZZ_{\geq0}^s$, there is $\alpha\in\calS$ with $u\leq\alpha$ if
and only if
\begin{equation}\label{eq:rank-extension}
 u(T)\leq r(T)\quad\text{for every }T\subseteq[s].
\end{equation}

We also use the numerical Hodge index theorem: for an $\RR$-Cartier divisor
$D$ and an ample $A$ on a normal projective variety of dimension at least two,
\begin{equation}\label{eq:hodge}
 DA^{n-1}=D^2A^{n-2}=0
\end{equation}
implies $D\equiv0$; see, for example, \cite[Lemma~3.2]{ZhangHodge}.

\section{The linear programme and finite-difference vanishing}\label{sec:optimal}

Throughout this section, $X$ is a complex klt Fano variety,
$n:=\dim X$, $\rho:=\rho(X)$, and
$\Lambda_X:=\Pic(X)/\!\equiv_{\mathrm{num}}$; $\tau_X$ and $c_X$ are as in
Definition~\ref{def:nef-complexity}.  Let
$\operatorname{Hilb}(\Nef(X)\cap\Lambda_X)=\{H_1,\ldots,H_N\}$ and choose nef
Cartier representatives.  Define $B:\RR^N\to N^1(X)_{\RR}$ by
$B(x)=\sum_{\alpha=1}^N x_\alpha[H_\alpha]$ and put
$\mathcal F:=B^{-1}([-K_X])\cap\RR_{\geq0}^N$.  We consider the linear programme
\begin{equation}\label{eq:LP}
 \text{maximize }\ell(x):=\sum_{\alpha=1}^N x_\alpha
 \quad\text{subject to }x\in\mathcal F.
\end{equation}

\begin{proposition}\label{prop:LP}
The maximum in \eqref{eq:LP} equals $\tau_X$.  It is attained at a rational
vertex, and the nonzero coordinates of such a vertex give a decomposition
\begin{equation}\label{eq:optimal}
 -K_X\equiv\sum_{i=1}^{s}a_iL_i,
 \qquad a_i\in\QQ_{>0},
 \qquad \sum_i a_i=\tau_X,
\end{equation}
where the numerical classes $[L_1],\ldots,[L_s]$ are linearly independent.
In particular, $s\leq\rho$.
\end{proposition}
\begin{proof}
The set $\mathcal F$ is nonempty because a sufficiently divisible multiple of
$-K_X$ is an ample Cartier divisor and hence has a Hilbert-basis expansion.
It is bounded: if $A$ is ample and
$\mu:=\min_\alpha H_\alpha A^{n-1}>0$, then every $x\in\mathcal F$ satisfies
\begin{equation}\label{eq:LP-bound}
 \mu\sum_{\alpha=1}^N x_\alpha
 \leq (-K_X)A^{n-1}.
\end{equation}
Thus $\mathcal F$ is a nonempty rational polytope, and
Lemma~\ref{lem:polyhedral-basic} gives a rational maximizing vertex.

Refining any admissible decomposition into Hilbert-basis elements
does not decrease its coefficient sum. Thus
$\tau_X=\max_{x\in\mathcal F}\ell(x)$. By
Lemma~\ref{lem:polyhedral-basic}, this maximum is attained at a
rational vertex, whose positive support has linearly independent
numerical classes.
\end{proof}

\begin{definition}[Optimal vertex decomposition of $-K_X$]\label{def:optimal-vertex-decomposition}
Let $v\in\mathcal F$ be a maximizing vertex in \eqref{eq:LP}.  The
decomposition \eqref{eq:optimal} determined by its positive coordinates is
called an \emph{optimal vertex decomposition} of $-K_X$.
\end{definition}

\begin{definition}[Multivariate Euler-characteristic polynomial]\label{def:multivariate-euler}
Let $X$ be a projective klt variety of dimension $n$, and let
$L_1,\ldots,L_s$ be Cartier divisors on $X$. Since klt singularities are
rational \cite[Theorem~5.22]{KM}, Hirzebruch--Riemann--Roch on a resolution
shows that there is a unique polynomial
\(P(t_1,\ldots,t_s)\in\QQ[t_1,\ldots,t_s]\)
of total degree at most $n$ such that
\[
 P(k)
 =\chi\left(X,\OO_X\left(-\sum_i k_iL_i\right)\right)
 \qquad(k\in\ZZ^s).
\]
We call $P$ the \emph{multivariate
Euler-characteristic polynomial} associated with $(X;L_1,\ldots,L_s)$.
\end{definition}

Fix an optimal vertex decomposition
\begin{equation}\label{eq:standing-optimal}
 -K_X\equiv\sum_{i=1}^{s}a_iL_i,
 \qquad a_i\in\QQ_{>0},
 \qquad \sum_i a_i=\tau_X,
\end{equation}
with linearly independent classes.  Define
\begin{equation}\label{eq:rounding}
 m_i:=\lceil a_i\rceil-1,
 \qquad
 \eps_i:=m_i+1-a_i,
 \qquad
 M:=\sum_i m_i,
 \qquad
 e:=\sum_i\eps_i.
\end{equation}
Then
\begin{equation}\label{eq:standing-rounded}
 -K_X\equiv\sum_{i=1}^{s}(m_i+1-\eps_i)L_i,
 \qquad
 m_i\in\ZZ_{\geq0},\quad 0\leq\eps_i<1.
\end{equation}
For $k=(k_1,\ldots,k_s)\in\ZZ^s$, put
\(D(k):=\sum_i k_iL_i.\)
Let $P$ denote the multivariate Euler-characteristic polynomial
associated with $(X;L_1,\ldots,L_s)$.

For a multi-index $u=(u_1,\ldots,u_s)$ we write
$t^u=t_1^{u_1}\cdots t_s^{u_s}$, and
$\operatorname{coeff}_{t^u}(P)$ denotes the coefficient of the monomial
$t^u$ in $P$.

For a polynomial $Q(t)$, let $T_i$ denote the translation operator in the
$i$-th coordinate, defined by
\((T_iQ)(t):=Q(t+\ee_i).\)
We use this notation throughout.

\begin{lemma}[Finite differences]\label{lem:finite-differences}
Let $P=P(t_1,\ldots,t_s)\in\QQ[t_1,\ldots,t_s]$.  Put
$\Delta_iP(t):=P(t+\ee_i)-P(t)$ and, for
$u=(u_1,\ldots,u_s)\in\ZZ_{\geq0}^s$, set
$\Delta^u:=\Delta_1^{u_1}\cdots\Delta_s^{u_s}$,
$|u|:=\sum_i u_i$, and $u!:=\prod_i u_i!$.  Then:
\begin{enumerate}[label=(\roman*)]
\item For every $t\in\ZZ^s$,
\begin{equation}\label{eq:finite-difference-formula}
 \Delta^uP(t)
 =\sum_{0\leq k\leq u}
 (-1)^{|u|-|k|}\binom{u}{k}P(t+k),
 \qquad
 \binom{u}{k}:=\prod_i\binom{u_i}{k_i}.
\end{equation}
\item If $\deg P<|u|$, then $\Delta^uP=0$.
\item If $\deg P\leq|u|$, then $\Delta^uP$ is constant and
\begin{equation}\label{eq:finite-difference-coefficient}
 \Delta^uP
 =u!\,\operatorname{coeff}_{t^u}(P).
\end{equation}
Equivalently, if $P_{|u|}$ is the homogeneous component of degree $|u|$,
then
\[
 \Delta^uP
 =u!\,\operatorname{coeff}_{t^u}(P_{|u|}).
\]
\end{enumerate}
\end{lemma}
\begin{proof}
Since $\Delta_i=T_i-\operatorname{id}$ and the $T_i$ commute, the
binomial theorem gives
\[
 \Delta^u
 =\sum_{0\leq k\leq u}
 (-1)^{|u|-|k|}\binom{u}{k}
 T_1^{k_1}\cdots T_s^{k_s},
\]
which is \eqref{eq:finite-difference-formula}.  Each $\Delta_i$ lowers the
total degree of a nonconstant polynomial by at least one, proving
(ii) and showing that $\Delta^uP$ is constant when
$\deg P\leq|u|$.  In total degree $|u|$, every monomial $t^\alpha$ with
$\alpha\neq u$ is killed by some $\Delta_i^{u_i}$, whereas
\[
 \Delta^ut^u=\prod_i u_i!=u!.
\]
This proves (iii).
\end{proof}

\begin{remark}\label{rem:finite-diff}
For integral $k$, $\Delta_i$ compares the Euler characteristics for $D(k)$
and $D(k)+L_i$.  If $L_i$ is globally generated and $H_i\in|L_i|$ is general,
the restriction sequence gives
\[
 \Delta_iP(k)=-\chi\bigl(H_i,\OO_{H_i}(-D(k)|_{H_i})\bigr).
\]
Thus, up to sign, $\Delta_i$ behaves like a virtual cut in the $L_i$-direction,
and iterated differences like successive cuts.  This is only geometric
intuition; the algebraic definition does not require $|L_i|$ to be base-point
free.
\end{remark}

\begin{lemma}[Vanishing box]\label{lem:Bm-vanishing}
Let $(X,\Delta)$ be a projective klt log Fano pair of dimension $n$ with
$-(K_X+\Delta)\equiv\sum_{i=1}^{s}a_iL_i$, where $a_i>0$ and the $L_i$ are
numerically nonzero nef Cartier divisors.  Put
$m_i:=\lceil a_i\rceil-1$, $M:=\sum_i m_i$, and
$D(k):=\sum_i k_iL_i$ for $k\in\ZZ^s$, and let $P$ be the associated
multivariate Euler-characteristic polynomial, so that
$P(k)=\chi(X,\OO_X(-D(k)))$.  Define the vanishing box
\[
 \mathcal B_m
 :=
 \left\{
 k\in\ZZ^s
 \ \middle|\
 0\leq k_i\leq m_i\text{ for all }i
 \right\}.
\]
Then
\begin{equation}\label{eq:Bm-vanishing}
 P(k)=
 \begin{cases}
  1,&k=0,\\
  0,&k\neq0,
 \end{cases}
 \qquad
 (k\in\mathcal B_m).
\end{equation}
Moreover,
\begin{equation}\label{eq:M-bound}
 M\leq n.
\end{equation}
\end{lemma}

\begin{proof}
For $k\in\mathcal B_m$,
\[
 -(K_X+\Delta)-D(k)
 \equiv\sum_i(a_i-k_i)L_i
\]
is ample.  Kawamata--Viehweg vanishing therefore gives
$H^q(X,\OO_X(-D(k)))=0$ for $q>0$.  If $k\neq0$, then $D(k)$ is a nonzero
nef class.  A nonzero section of $\OO_X(-D(k))$ would give a nonzero effective
divisor numerically equivalent to $-D(k)$, contradicting its intersection
with an ample $(n-1)$-fold product.  Thus \eqref{eq:Bm-vanishing} holds.

Apply Lemma~\ref{lem:finite-differences} (i) with $u=m$.  Only the
term $k=0$ survives, so
\begin{equation}\label{eq:Bm-difference}
 \Delta^mP(0)=(-1)^M\neq0.
\end{equation}
If $M>n$, then $\deg P\leq n<M$, contradicting
Lemma~\ref{lem:finite-differences} (ii).  Hence $M\leq n$.
\end{proof}

\begin{proposition}\label{prop:defect}
For the fixed optimal vertex decomposition,
\begin{equation}\label{eq:defect}
 c_X=(n-M)+(\rho-s)+e\geq0.
\end{equation}
In particular, if $c_X<1$, then
\[
 M=n,\qquad s=\rho,\qquad e=c_X<1.
\]
\end{proposition}
\begin{proof}
By optimality and \eqref{eq:rounding},
$\tau_X=M+s-e$.  The equality in \eqref{eq:defect} follows immediately,
and its three summands are nonnegative by
Lemma~\ref{lem:Bm-vanishing} and $s\leq\rho$.
\end{proof}

\section{Finite differences and the equality case}\label{sec:equality}

Throughout this section, $X$ is a complex klt Fano variety of dimension
$n$, and we retain the optimal vertex decomposition and notation from
Section~\ref{sec:optimal}:
\[
 -K_X\equiv\sum_{i=1}^{s}(m_i+1-\eps_i)L_i,
 \quad m=(m_1,\ldots,m_s),\quad M=\sum_i m_i,\quad e=\sum_i\eps_i.
\]
For $t=(t_1,\ldots,t_s)$, write $D(t)=\sum_i t_iL_i$ and
$L^\alpha=\prod_iL_i^{\alpha_i}$, let $P(t)$ be the associated multivariate
Euler-characteristic polynomial, and put
$\calS=\{\alpha\in\ZZ_{\geq0}^s\mid|\alpha|=n,\ L^\alpha>0\}$ and
$\mathcal B_m=\{k\in\ZZ_{\geq0}^s\mid0\leq k_i\leq m_i\}$.  We also retain
$\Delta_iP(t)=P(t+\ee_i)-P(t)$ and
$\Delta^u=\Delta_1^{u_1}\cdots\Delta_s^{u_s}$.

A normal Fano curve is $\PP^1$, with $\tau=2$ and $c=0$.
We henceforth treat $n\geq2$; the elementary top-difference statements
below also apply in dimension one.

\subsection{The two leading terms of Riemann--Roch}

Since klt singularities are rational \cite[Theorem~5.22]{KM}, a resolution
$\pi:\widetilde X\to X$ gives
\[
 \chi\bigl(X,\OO_X(-D)\bigr)
 =
 \chi\bigl(\widetilde X,\OO_{\widetilde X}(-\pi^*D)\bigr)
\]
for every Cartier divisor $D$.

Write
\(P(t)=P_n(t)+P_{n-1}(t)+P_{\leq n-2}(t),\)
where $P_j(t)$ denotes the homogeneous component of $P(t)$ of total
degree $j$, and $P_{\leq n-2}(t)$ denotes the sum of its homogeneous
components of total degree at most $n-2$.

Hirzebruch--Riemann--Roch on $\widetilde X$ gives the two leading terms
\[
 P_n(t)=(-1)^n\frac{D(t)^n}{n!},\qquad
 P_{n-1}(t)=(-1)^n\frac{K_XD(t)^{n-1}}{2(n-1)!}.
\]
Hence
\begin{equation}\label{eq:RR}
 P(t)=(-1)^n\left(
 \frac{D(t)^n}{n!}
 +\frac{K_XD(t)^{n-1}}{2(n-1)!}
 \right)+P_{\leq n-2}(t).
\end{equation}
Indeed, the degree-one part of the Todd class is $-K_{\widetilde X}/2$,
and exceptional divisors have zero intersection with $n-1$ pullbacks of
Cartier divisors from $X$.

For a polynomial $Q(t)$, write
$\partial_i=\partial/\partial t_i$ and $\partial^u=\prod_i\partial_i^{u_i}$.
The translation operator is $T_i=e^{\partial_i}$ on polynomials, hence
$\Delta_i=e^{\partial_i}-1$.  Consequently,
\[
 (e^z-1)^r=z^r+\frac r2z^{r+1}+O(z^{r+2})
\]
gives, for $|u|=n-1$,
\begin{equation}\label{eq:difference-operator-expansion}
 \Delta^u
 =\partial^u+\frac12\sum_i u_i\partial^{u+\ee_i}
 +\text{operators of order at least }n+1.
\end{equation}
To evaluate \eqref{eq:difference-operator-expansion} on $P(t)$ at the
origin, first note that, for every multi-index $\alpha$ with
$|\alpha|=m$,
\[
 \left.
 \partial^\alpha\!\left(\frac{D(t)^m}{m!}\right)
 \right|_{t=0}
 =L^\alpha.
\]
Indeed, the multinomial expansion gives
\[
 D(t)^m
 =
 \sum_{|\beta|=m}\frac{m!}{\beta!}t^\beta L^\beta,
\]
so after applying $\partial^\alpha$ and evaluating at $t=0$, only the
term with $\beta=\alpha$ remains.

Now suppose $|u|=n-1$. Since $\partial^u$ has order $n-1$, it
annihilates $P_{\leq n-2}(t)$. Its application to $P_n(t)$ is
homogeneous of degree one and therefore vanishes at $t=0$, whereas
\[
 \partial^uP_{n-1}(0)
 =
 \frac{(-1)^n}{2}K_XL^u.
\]
On the other hand, each $\partial^{u+\ee_i}$ has order $n$, so it
annihilates both $P_{n-1}(t)$ and $P_{\leq n-2}(t)$, while
\[
 \partial^{u+\ee_i}P_n(0)
 =
 (-1)^nL^{u+\ee_i}
 =
 (-1)^nL_iL^u.
\]
Finally, every operator of order at least $n+1$ annihilates $P(t)$,
since $P(t)$ has total degree at most $n$. Hence
\[
 \Delta^uP(0)
 =
 \frac{(-1)^n}{2}K_XL^u
 +
 \frac{(-1)^n}{2}
 \left(\sum_i u_iL_i\right)L^u.
\]
Therefore
\begin{equation}\label{eq:RR-difference}
 2(-1)^n\Delta^uP(0)
 =\left(K_X+\sum_i u_iL_i\right)L^u
 \qquad(|u|=n-1).
\end{equation}

Since finite differences are integral combinations of Euler
characteristics, \eqref{eq:RR-difference} also gives
\begin{equation}\label{eq:parity}
 \left(K_X+\sum_i u_iL_i\right)L^u\in2\ZZ
 \qquad(|u|=n-1).
\end{equation}
The remaining summand is an integral Cartier intersection; hence
\begin{equation}\label{eq:K-integral}
 K_XL^u\in\ZZ\qquad(|u|=n-1).
\end{equation}

\begin{remark}
The calculation above is a multivariate version of the standard
one-variable relation between the Hilbert polynomial, finite differences,
and adjunction.  When $s=1$, writing $L=L_1$ and $u=n-1$,
\eqref{eq:RR-difference} becomes
\[
 2(-1)^n\Delta^{n-1}P(0)
 =
 \bigl(K_X+(n-1)L\bigr)L^{n-1}.
\]
In the smooth polarized case, the right-hand side is the usual adjunction
expression $2g(X,L)-2$ for the sectional genus.  For several nef divisors,
the same calculation replaces $(n-1)L$ by $\sum_i u_iL_i$ and
$L^{n-1}$ by the mixed intersection $L^u$, so that the finite difference
$\Delta^uP(0)$ detects the mixed adjunction intersection
\[
 \left(K_X+\sum_i u_iL_i\right)L^u.
\]
This viewpoint also explains the use of finite differences rather than
formal differentiation alone: the degree-$n$ and degree-$(n-1)$ parts of
Riemann--Roch determine the expression in \eqref{eq:RR-difference}, while
finite differences are integral combinations of Euler characteristics.
The latter property yields the parity and integrality statements
\eqref{eq:parity} and \eqref{eq:K-integral}.
\end{remark}

The next two lemmas are not restricted to the fixed optimal decomposition.
More generally, let $(X,\Delta)$ be a projective klt log Fano pair and let
\[
 -(K_X+\Delta)\equiv\sum_{i=1}^{s}a_iL_i,
 \qquad a_i>0,
\]
be any positive decomposition by nonzero nef Cartier divisors. Set
\[
 m_i=\lceil a_i\rceil-1,\qquad
 \eps_i=m_i+1-a_i,\qquad
 M=\sum_i m_i,\qquad e=\sum_i\eps_i,
\]
and define $P$, $L^\alpha$, and $\calS$ from these $L_i$ exactly as in
Section~\ref{sec:optimal}. Lemma~\ref{lem:Bm-vanishing} and the
finite-difference identities above apply verbatim to this decomposition.
After these standalone statements we return to the fixed optimal
decomposition.

\begin{lemma}\label{lem:full-degree}
Let $(X,\Delta)$ be a projective klt log Fano pair of dimension $n$, and suppose
\[
 -(K_X+\Delta)\equiv\sum_{i=1}^{s}a_iL_i,
 \qquad a_i>0,
\]
where the $L_i$ are numerically nonzero nef Cartier divisors. Define
\[
 m_i:=\lceil a_i\rceil-1,
 \qquad
 \eps_i:=m_i+1-a_i,
 \qquad
 M:=\sum_i m_i,
 \qquad
 e:=\sum_i\eps_i,
\]
and write $L^\alpha:=\prod_iL_i^{\alpha_i}$. If $M=n$, then
$L^m=1$. For every $j$ with $m_j>0$, set
$C_j=L^{m-\ee_j}$ and $I_{ij}=L_iC_j$. Then
\begin{equation}\label{eq:epsilon}
 \eps_j=\sum_{i\neq j}(1-\eps_i)I_{ij}+\Delta C_j.
\end{equation}
If $e<1$, then
\[
 I_{ij}=L_iC_j=L^{m-\ee_j+\ee_i}=0
 \qquad
 \text{for every $j$ with $m_j>0$ and every $i\ne j$},
\]
and
\[
 \eps_j=\Delta C_j
 \qquad\text{for every $j$ with $m_j>0$}.
\]
In particular, when $\Delta=0$, equation \eqref{eq:epsilon} becomes
\[
 \eps_j=\sum_{i\neq j}(1-\eps_i)I_{ij},
\]
and $e<1$ implies $\eps_j=0$ for every $j$ with $m_j>0$.
\end{lemma}
\begin{proof}
Since $M=n$, one has $|m|=n$. By
Lemma~\ref{lem:finite-differences} (iii) and the degree-$n$
homogeneous part of \eqref{eq:RR},
\begin{align*}
 \Delta^mP(0)
 &=m!\,\operatorname{coeff}_{t^m}\!\left(
   \frac{(-1)^n}{n!}\left(\sum_i t_iL_i\right)^n\right)\\
 &=(-1)^nL^m.
\end{align*}
Indeed,
\[
 \left(\sum_i t_iL_i\right)^n
 =\sum_{|\alpha|=n}\frac{n!}{\alpha!}t^\alpha L^\alpha,
\]
so the coefficient of $t^m$ in the degree-$n$ part of $P$ is
$(-1)^nL^m/m!$. On the other hand, \eqref{eq:Bm-difference} gives
\[
 \Delta^mP(0)=(-1)^M=(-1)^n.
\]
Comparing the two expressions yields
\(L^m=1.\)

Fix $j$ with $m_j>0$ and put $u=m-\ee_j$. Then $|u|=n-1$.
By Lemma~\ref{lem:finite-differences} (i),
\[
 \Delta^uP(0)
 =\sum_{0\leq k\leq u}
   (-1)^{|u|-|k|}\binom{u}{k}P(k).
\]
Every $k$ occurring in this sum lies in $\mathcal B_m$. Hence
Lemma~\ref{lem:Bm-vanishing} shows that only $k=0$ remains, and therefore
\[
 \Delta^uP(0)=(-1)^{n-1}.
\]
The calculation leading to \eqref{eq:RR-difference} remains valid for the
underlying klt variety even if $K_X$ is not separately $\QQ$-Cartier: on a
resolution the degree-$(n-1)$ term pushes forward to the canonical Weil
cycle on $X$, whose intersections with Cartier divisors are defined. Thus
\[
 \left(K_X+\sum_i(m_i-\delta_{ij})L_i\right)C_j=-2.
\]
Since
\[
 K_X+\Delta\equiv-\sum_i(m_i+1-\eps_i)L_i,
\]
we obtain
\[
 \left((K_X+\Delta)+\sum_i(m_i-\delta_{ij})L_i\right)C_j
 =-(2-\eps_j)L_jC_j
  -\sum_{i\ne j}(1-\eps_i)L_iC_j.
\]
Using $L_jC_j=L^m=1$ and subtracting $\Delta C_j$ gives
\eqref{eq:epsilon}.

Every term on the right of \eqref{eq:epsilon} is nonnegative. If
$I_{ij}>0$ for some $i\ne j$, then $I_{ij}\geq1$, so
$\eps_j\geq1-\eps_i$ and hence $\eps_i+\eps_j\geq1$, contradicting
$e<1$. Thus all such $I_{ij}$ vanish, and \eqref{eq:epsilon} gives
$\eps_j=\Delta C_j$. The final assertion is the special case
$\Delta=0$.
\end{proof}

\begin{lemma}\label{lem:singleton}
Let $(X,\Delta)$ be a projective klt log Fano pair of dimension $n$, and suppose
\[
 -(K_X+\Delta)\equiv\sum_{i=1}^{s}a_iL_i,
 \qquad a_i>0,
\]
where the $L_i$ are numerically nonzero nef Cartier divisors. Define
$m_i,\eps_i,M,e$ as in Lemma~\ref{lem:full-degree}, write
$L^\alpha:=\prod_iL_i^{\alpha_i}$, and set
\[
 \calS:=\left\{\alpha\in\ZZ_{\geq0}^s\ \middle|\ |\alpha|=n,\ L^\alpha>0\right\}.
\]
Assume $M=n$ and
\[
 I_{ij}=L^{m-\ee_j+\ee_i}=0
 \qquad
 \text{for every $j$ with $m_j>0$ and every $i\ne j$}.
\]
Then every $m_i>0$ and $\calS=\{m\}$.
\end{lemma}
\begin{proof}
By Lemma~\ref{lem:full-degree}, $L^m=1$, so $m\in\calS$.
Suppose first that $m_i=0$ for some $i$. Since $L_i$ is nonzero nef
and $A=\sum_rL_r$ is ample, one has $L_iA^{n-1}>0$. Expanding,
\[
 L_iA^{n-1}
 =\sum_{|\gamma|=n-1}\frac{(n-1)!}{\gamma!}L^{\gamma+\ee_i}.
\]
All terms are nonnegative, hence $L^{\gamma+\ee_i}>0$ for some
$\gamma$. Thus, with $\beta=\gamma+\ee_i$, we have
$\beta\in\calS$ and $\beta_i\geq1>0=m_i$. By symmetric exchange applied to
$\beta$ and $m$, there is an index $j$ with $\beta_j<m_j$ and
\(m-\ee_j+\ee_i\in\calS.\)
Here $j\ne i$ and $m_j>0$. Therefore
\(I_{ij}=L^{m-\ee_j+\ee_i}>0.\)
This contradicts the assumed equality $I_{ij}=0$, which applies because
$m_j>0$ and $i\ne j$. Hence $m_i>0$ for every $i$.

Now suppose that $\beta\in\calS$ and $\beta\ne m$. Since
$|\beta|=|m|=n$, there is an index $j$ with $m_j>\beta_j$.
Applying symmetric exchange to $m$ and $\beta$, we obtain an index
$i$ with $m_i<\beta_i$ and
\(m-\ee_j+\ee_i\in\calS.\)
Again $i\ne j$ and $m_j>0$, so
$0<I_{ij}=L^{m-\ee_j+\ee_i}$. Since $m_j>0$ and $i\ne j$,
this contradicts the assumed equality $I_{ij}=0$. Thus no such
$\beta$ exists and $\calS=\{m\}$.
\end{proof}

We recall the following degree-one consequence of Fujita's
classification of polarized varieties of $\Delta$-genus zero
\cite{Fujita75}. For a normal projective $d$-fold $Y$ and an ample
Cartier divisor $H$, put
$\Delta(Y,H)=d+H^d-h^0(Y,H)$. If
\begin{equation}\label{eq:Fujita-degree-one}
 H^d=1\quad\text{and}\quad h^0(Y,H)=d+1,
\end{equation}
then $(Y,H)\simeq(\PP^d,\OO(1))$.

\begin{lemma}\label{lem:image-P}
Let $(X,\Delta)$ be a projective klt log Fano pair,
and let $f:X\to Y$ be the contraction of a nonzero nef Cartier divisor
$L$. Write $L\sim f^*H$ with $H$ ample Cartier, and put $d=\dim Y$.
Suppose that either
\begin{enumerate}[label=(\roman*)]
\item $\chi(X,-tL)=0$ for $t=1,\ldots,d$; or
\item $H^d=1$ and $\chi(X,-tL)=0$ for $t=1,\ldots,d-1$.
\end{enumerate}
Then $(Y,H)\simeq(\PP^d,\OO(1))$.
\end{lemma}
\begin{proof}
By Lemma~\ref{lem:descent},
\[
 Q(t):=\chi(Y,\OO_Y(-tH))=\chi(X,\OO_X(-tL))
\]
is a polynomial of degree $d$. Kawamata--Viehweg vanishing for
$(X,\Delta)$ gives $H^q(X,\OO_X)=0$ for $q>0$, and hence $Q(0)=1$.
Its leading term is
\[
 Q(t)=\frac{(-1)^dH^d}{d!}t^d+O(t^{d-1}).
\]
Under (i), the roots $1,\ldots,d$ and the value $Q(0)=1$
give
\begin{equation}\label{eq:P-polynomial}
 Q(t)=\prod_{j=1}^{d}\left(1-\frac{t}{j}\right),
\end{equation}
and comparison of leading coefficients gives $H^d=1$.
Under (ii), let $R(t)$ denote the product on the right of
\eqref{eq:P-polynomial}. The leading coefficients of $Q$ and $R$
agree, so $\deg(Q-R)\leq d-1$. Moreover, $Q-R$ vanishes at the $d$
distinct integers $0,1,\ldots,d-1$. Hence $Q=R$, and
\eqref{eq:P-polynomial} holds in this case as well. This also covers
$d=1$, when the root condition in (ii) is empty.

Moreover,
\[
 L-(K_X+\Delta)=L+\bigl(-(K_X+\Delta)\bigr)
\]
is ample. Kawamata--Viehweg vanishing therefore yields
$H^q(X,\OO_X(L))=0$ for $q>0$. Since $L\sim f^*H$ and
$f_*\OO_X=\OO_Y$, the projection formula gives
\[
 H^0(X,\OO_X(L))\simeq H^0(Y,\OO_Y(H)).
\]
Therefore
\[
 h^0(Y,H)=\chi(X,L)=\chi(Y,H)=Q(-1)
 =\prod_{j=1}^{d}\left(1+\frac1j\right)=d+1.
\]
Thus $H^d=1$ and $h^0(Y,H)=d+1$. By the degree-one case above,
$(Y,H)\simeq(\PP^d,\OO_{\PP^d}(1))$.
\end{proof}

\begin{proposition}\label{prop:product}
Let $(X,\Delta)$ be a projective klt log Fano pair of dimension $n>0$.
Suppose that
\[
 -(K_X+\Delta)\equiv
 \sum_{i=1}^{s}(m_i+1-\eps_i)L_i,
 \qquad
 m_i\in\ZZ_{\geq0},\quad0\leq\eps_i<1,
\]
where the $L_i$ are nonzero nef Cartier divisors and
\[
 \sum_i m_i=n,
 \qquad
 \sum_i\eps_i<1.
\]
Then every $m_i>0$ and
\[
 (X;L_1,\ldots,L_s)\simeq
 \left(\prod_{i=1}^{s}\PP^{m_i};\ p_1^*\OO(1),\ldots,p_s^*\OO(1)\right).
\]
Under this identification,
\begin{equation}\label{eq:product-boundary}
 \Delta\equiv\sum_i\eps_iL_i.
\end{equation}
In particular, $\rho(X)=s$.
\end{proposition}
\begin{proof}
Lemma~\ref{lem:full-degree} gives $L^m=1$ and
\[
 L^{m-\ee_j+\ee_i}=0
 \qquad
 \text{for every $j$ with $m_j>0$ and every $i\ne j$}.
\]
Lemma~\ref{lem:singleton} therefore gives $m_i>0$ for every $i$ and
$\calS=\{m\}$. Formula \eqref{eq:numdim} gives
$\nu(L_i)=m_i$. The contraction $f_i$ of $L_i$ has image of
dimension $m_i$; Lemma~\ref{lem:Bm-vanishing} and Lemma~\ref{lem:image-P}
identify it with $\PP^{m_i}$ and give $L_i=f_i^*\OO(1)$.

Let
\[
 f=(f_1,\ldots,f_s):X\longrightarrow
 Y:=\prod_{i=1}^{s}\PP^{m_i},
\]
and let $H_i=p_i^*\OO_{\PP^{m_i}}(1)$. Then $L_i=f^*H_i$ for every
$i$, and
\[
 f^*\OO_Y(1,\ldots,1)=\sum_iL_i
\]
is ample. Hence $f$ has no positive-dimensional fibre: otherwise a
curve $C$ contained in such a fibre would satisfy
$(\sum_iL_i)\cdot C=0$, contradicting ampleness. Thus $f$ is
quasi-finite, and since it is projective, it is finite.

Moreover,
\[
 \dim Y=\sum_i m_i=n=\dim X,
\]
so the finite image $f(X)$ is all of $Y$. By the projection formula,
\[
 L_1^{m_1}\cdots L_s^{m_s}
 =\deg(f)\,H_1^{m_1}\cdots H_s^{m_s}
 =\deg(f),
\]
because
\(H_1^{m_1}\cdots H_s^{m_s}=1\)
on $\prod_i\PP^{m_i}$. Since Lemma~\ref{lem:full-degree} gives
$L_1^{m_1}\cdots L_s^{m_s}=1$, we obtain $\deg(f)=1$.
Thus $f$ is finite and birational. Since $Y$ is normal, $f$ is an
isomorphism.

Under this identification,
\[
 K_X\equiv-\sum_i(m_i+1)L_i,
\]
whereas the assumed decomposition gives
\[
 K_X+\Delta\equiv-\sum_i(m_i+1-\eps_i)L_i.
\]
Subtracting these identities yields \eqref{eq:product-boundary}.
\end{proof}

\begin{theorem}\label{thm:gap}
Let $X$ be a complex klt Fano variety. Then
$c_X<1$ if and only if $X$ is a product of projective spaces. In this
case $c_X=0$.
\end{theorem}
\begin{proof}
If $c_X<1$, Proposition~\ref{prop:defect} gives $M=n$ and $e<1$.
Applying Proposition~\ref{prop:product} to the optimal decomposition with
$\Delta=0$ shows that $X$ is a product of projective spaces.
Conversely, on $X=\prod_i\PP^{d_i}$ the standard decomposition
$-K_X=\sum_i(d_i+1)p_i^*\OO(1)$ has coefficient sum $n+\rho$.
Proposition~\ref{prop:defect} gives the reverse inequality and
therefore $c_X=0$.
\end{proof}

\begin{lemma}\label{lem:index}
Let $X$ be a klt Fano variety of dimension $n$ and let $H$ be an ample
Cartier divisor. If
\[
 -K_X\equiv tH
\]
for a positive integer $t$, then $t\leq n+1$. If $t=n+1$, then
$(X,H)\simeq(\PP^n,\OO(1))$.
\end{lemma}
\begin{proof}
If $t\leq n$ there is nothing to prove. If $t>n$, numerical and
$\QQ$-linear equivalence agree on a klt Fano variety by
Lemma~\ref{lem:descent}. Hence Fujino--Miyamoto
\cite[Theorem~1.1]{FujinoMiyamoto} applies and gives
$(X,H)\simeq(\PP^n,\OO(1))$. Consequently $t=n+1$.
\end{proof}

\begin{remark}\label{rem:singular-KO}
Araujo--Druel \cite[Theorem~1.1]{AraujoDruel} prove the corresponding
singular Kobayashi--Ochiai statement in greater generality. In particular,
for a $\QQ$-Fano variety the next value $-K_X\sim_{\QQ}nH$ gives a possibly
singular quadric.
\end{remark}

\section{Full-degree decompositions with small fractional defect}\label{frac:sec:full}

Here $X$ is a complex klt Fano variety of dimension
$n$.  Unless a statement explicitly introduces local notation, we retain
from Sections~\ref{sec:optimal}--\ref{sec:equality} the standing optimal
vertex decomposition, the rounding data $m_i,\eps_i,M,e$, the vector
$m=(m_1,\ldots,m_s)$, and the notation
\[
 L^\alpha:=\prod_iL_i^{\alpha_i},
 \qquad
 \calS=\{\alpha\in\ZZ_{\geq0}^s\mid |\alpha|=n,\ L^\alpha>0\},
 \qquad
 \mathcal B_m=\{k\in\ZZ_{\geq0}^s\mid 0\leq k_i\leq m_i\}.
\]

\subsection{A birational lemma over a product}\label{frac:sec:regularity}

\begin{lemma}\label{frac:lem:star}
Let $X$ be a complex klt Fano variety of dimension $n$, and let
$f:X\to Y=\prod_{i=1}^{q}\PP^{m_i}$ be a birational morphism, where
$m_i>0$ and $\sum_i m_i=n$. For each $i$, let $H_i$ be the pullback
to $Y$ of a hyperplane divisor on the $i$-th factor, and put
$L_i=f^*H_i$. Let $L$ be a nef Cartier divisor on $X$, and suppose
\[
 -K_X\equiv bL+\sum_i(m_i+1-bh_i)L_i,
 \qquad0<b<1,\quad h_i\in\ZZ_{\geq0},\quad bh_i<1.
\]
Then $f$ is an isomorphism and
\[
 L\sim\sum_{i=1}^q h_iL_i.
\]
In particular, $L,L_1,\ldots,L_q$ cannot be numerically independent.
\end{lemma}
\begin{proof}
Put $H=\sum_i h_iH_i$ and $\lambda=1/b$, and let
$\{E_\alpha\}$ be all the $f$-exceptional prime divisors on $X$.
Since $Y$ is smooth, the relative canonical divisor is
$R:=K_X-f^*K_Y=\sum_\alpha a_\alpha E_\alpha$ with
$a_\alpha\in\ZZ_{>0}$.
The assumed identity gives
\begin{equation}\label{frac:eq:star-relative-canonical}
 R\equiv b(f^*H-L).
\end{equation}
Choose a Cartier divisor representing the class $L$, still denoted by
$L$, and put $G=f_*L$. This is a Weil divisor on the smooth variety
$Y$, hence Cartier, and $D:=L-f^*G$ is an exceptional Cartier divisor.
The divisor $D+\lambda R$ is exceptional and numerically trivial over
$Y$. Applying the negativity lemma
\cite[Lemma~4.16]{FujinoFundamental} to it and to its negative gives
$D+\lambda R=0$. Thus $L=f^*G-\lambda R$, and
\eqref{frac:eq:star-relative-canonical} implies $f^*(G-H)\equiv0$.
The injectivity of $f^*$ on numerical divisor classes gives $G\equiv H$.
Since numerical and linear equivalence agree on $Y$, we obtain
\begin{equation}\label{frac:eq:star-divisor}
 L\sim f^*H-F,\qquad F:=\lambda R=-D,
\end{equation}
where $F$ is effective and, by $F=-D$, integral and Cartier.

Suppose that $f$ has an exceptional prime, and put
$Z=\bigcup_\alpha f(E_\alpha)$. Purity of the exceptional locus
\cite[Chapter~2, Corollary~2.63]{KM} applies because $X$ is projective
and $Y$ is smooth, hence $\QQ$-factorial. It gives
$\operatorname{Exc}(f)=\bigcup_\alpha E_\alpha=\Supp R$.
Thus $Z=f(\operatorname{Exc}(f))$ is precisely the locus over which
$f$ is not an isomorphism, and it is nonempty of codimension at least
two. In particular $n\geq2$. If $H=0$, then $L\sim-F$ contradicts nefness, since $F$ is
nonzero and effective.
Thus $h_{\max}:=\max_i h_i>0$.

By Lemma~\ref{lem:nef-monoid}, choose a positive integer $r$ such that
$|rL|$ is base-point free. Pushing forward its divisors gives a linear
system $\mathcal M_r\subset|rH|=|\OO_Y(rh_1,\ldots,rh_q)|$.
For every $M\in\mathcal M_r$, we have $f^*M=B_M+rF$, where
$B_M$ runs through $|rL|$.
For every exceptional prime divisor $E_\alpha$, define
\[
 \ord_{E_\alpha}(\mathcal M_r)
 :=\min_{M\in\mathcal M_r}\operatorname{coeff}_{E_\alpha}(f^*M).
\]
Since $f^*M=B_M+rF$ and $|rL|$ is base-point free, some $B_M$
does not contain $E_\alpha$. Thus the minimum coefficient of
$E_\alpha$ in $B_M$ is zero, and
\begin{equation}\label{frac:eq:star-valuation}
 \ord_{E_\alpha}(\mathcal M_r)=r\lambda a_\alpha.
\end{equation}
A positive order along $E_\alpha$ means that every member vanishes
along its centre $f(E_\alpha)$. Conversely, $f$ is an isomorphism
over $Y\setminus Z$, where the base-point-freeness of $|rL|$ gives
no base points downstairs. Thus $\operatorname{Bs}(\mathcal M_r)=Z$
set-theoretically; in particular, $\mathcal M_r$ is movable.

We recall the convention for singularities of a movable linear
system \cite[Definition~8]{KollarRigidity}. For such a system $\mathcal M$
on a smooth variety $V$ and $c\geq0$, the expression
$(V,c\mathcal M)$ is a formal pair. For an exceptional prime divisor
$E$ over $V$, its discrepancy is
\[
 a(E;V,c\mathcal M):=a(E;V)-c\ord_E(\mathcal M),
\]
where $a(E;V)$ is its usual discrepancy over $V$ and
$\ord_E(\mathcal M)$ is the coefficient of $E$ in the fixed part of
the pullback system on a log resolution. The pair is \emph{canonical}
if these discrepancies are nonnegative for every exceptional divisor
over $V$. Moreover,
$\operatorname{mult}_y\mathcal M:=\min_{M\in\mathcal M}
\operatorname{mult}_y M$ is the order of the base ideal of
$\mathcal M$ at $y$.
By \cite[\S~13.7, P3]{BEV}, the function
$y\mapsto\operatorname{mult}_y\mathcal M_r$ is upper semicontinuous.

Since $\lambda>h_{\max}$, \eqref{frac:eq:star-valuation} shows that the
discrepancy of each $E_\alpha$ for
$(Y,(rh_{\max})^{-1}\mathcal M_r)$ satisfies
\[
 a_\alpha-\frac{\ord_{E_\alpha}(\mathcal M_r)}{rh_{\max}}
 = a_\alpha\left(1-\frac{\lambda}{h_{\max}}\right)<0.
\]
If $\operatorname{mult}_y\mathcal M_r\leq rh_{\max}$ at some
$y\in Z$, upper semicontinuity of multiplicity gives an open
neighbourhood $U$ of $y$ on which the same bound holds at every point.
The multiplicity estimate \cite[Lemma~28]{KollarRigidity}, applied on
$U$ with coefficient $(rh_{\max})^{-1}$, makes the restricted pair
canonical. This contradicts the negative discrepancy above for an
$E_\alpha$ whose centre contains $y$. Therefore, for every $y\in Z$,
\begin{equation}\label{frac:eq:star-corti}
 \operatorname{mult}_y\mathcal M_r>rh_{\max}\geq rh_i
 \qquad\text{for every }i.
\end{equation}
For $y=(y_1,\ldots,y_q)$, let $P_i(y)\simeq\PP^{m_i}$ be the
coordinate fibre obtained by fixing all coordinates except the $i$-th.
If a member $M\in\mathcal M_r$ does not contain $P_i(y)$, its
defining section restricts to a nonzero section of
$\OO_{P_i(y)}(rh_i)$. Restriction cannot decrease the vanishing order
at $y$, so this section vanishes to order at least
$\operatorname{mult}_y\mathcal M_r>rh_i$. This is impossible: on
an affine chart, a nonzero section of $\OO_{\PP^{m_i}}(rh_i)$ is
a polynomial of degree at most $rh_i$ and has vanishing order at
most $rh_i$. Thus every member contains $P_i(y)$, and hence
$P_i(y)\subset Z$ for every $y\in Z$ and every $i$. Starting from
any point of the nonempty set $Z$ and changing coordinates one at a
time reaches every point of $Y$. This gives $Z=Y$, a contradiction.

Therefore $f$ has no exceptional prime, and purity again makes it an
isomorphism. Now $R=0$, and \eqref{frac:eq:star-divisor} gives
$L\sim\sum_i h_iL_i$, proving also the asserted numerical dependence.
\end{proof}

\subsection{The full-degree argument}

\begin{proposition}\label{frac:prop:full}
Let $X$ be a complex klt Fano variety of dimension
$n$. Fix an optimal vertex decomposition \eqref{eq:standing-optimal}, with
the rounding data $m_i,\eps_i,M,e$ from \eqref{eq:rounding}. If
$M=n$ and $e<2$, then $e=0$ and $X$ is a product of projective spaces.
In particular $s=\rho(X)$ and $c_X=0$.
\end{proposition}
\begin{proof}
By Lemma~\ref{lem:full-degree},
\(L^m=1,\)
and, for every $j$ with $m_j>0$,
\begin{equation}\label{frac:eq:full-epsilon}
 \eps_j=\sum_{i\ne j}(1-\eps_i)I_{ij},
 \qquad
 I_{ij}:=L^{m-\ee_j+\ee_i}\in\ZZ_{\geq0}.
\end{equation}
Recall also that
\[
 \calS=\left\{\alpha\in\ZZ_{\geq0}^s\ \middle|\ |\alpha|=n,\ L^\alpha>0\right\}
\]
is $M$-convex and that $m\in\calS$ because $L^m=1$.

\smallskip\noindent
\emph{Step 1: the directed exchange graph.}
Let $G_m$ be the directed graph with vertex set
\([s]=\{1,\ldots,s\}.\)
Thus each index $i\in[s]$ is identified with the corresponding vertex
(or node) of $G_m$. For two distinct vertices $i,j\in[s]$, draw an arrow
\[
 j\longrightarrow i
\]
if
\(m_j>0\) and \(I_{ij}>0\).
Since $I_{ij}=L^{m-\ee_j+\ee_i}$, this is equivalent to
\[
 j\longrightarrow i
 \quad\text{if and only if}\quad
 m_j>0,\quad i\ne j,\quad m-\ee_j+\ee_i\in\calS.
\]
Thus an arrow $j\to i$ records an exchange in which one unit is removed
from the $j$-th coordinate of $m$ and added to its $i$-th coordinate. 
The out-degree of a vertex $j$ is the number of indices $i\ne j$
for which there is an arrow $j\to i$; equivalently,
\[
 \operatorname{outdeg}(j)
 =
 \#\{\,i\ne j\mid I_{ij}>0\,\}.
\]
We shall prove that $G_m$ has no arrows.

We first record four consequences of \eqref{frac:eq:full-epsilon} and
symmetric exchange.

Suppose that $j\to i$. Then $I_{ij}$ is a positive integer, hence
$I_{ij}\geq1$. Equation \eqref{frac:eq:full-epsilon}, applied to the
source vertex $j$, gives
\(\eps_j
 \geq (1-\eps_i)I_{ij}
 \geq 1-\eps_i.\)
Consequently
\begin{equation}\label{frac:eq:edge-defect}
 \eps_i+\eps_j\geq1.
\end{equation}
Both endpoint defects are in fact positive. Indeed,
$1-\eps_i>0$ gives $\eps_j>0$, while $\eps_i=0$ would imply
$\eps_j\geq1$, contrary to $\eps_j<1$.

Next, no vertex has two distinct outgoing arrows. Indeed, suppose
\(j\to i\) and \(j\to k\), with \(i\ne k\).
The two corresponding terms occur in the sum
\eqref{frac:eq:full-epsilon} for $\eps_j$, and therefore
\[
 \eps_j
 \geq (1-\eps_i)I_{ij}+(1-\eps_k)I_{kj}
 \geq (1-\eps_i)+(1-\eps_k).
\]
Hence
\(\eps_i+\eps_j+\eps_k\geq2,\)
and so $e=\sum_r\eps_r\geq2$, contradicting the assumption $e<2$.
Thus every vertex of $G_m$ has out-degree at most one.

Similarly, two arrows cannot have disjoint sets of endpoint vertices.
Indeed, if $j\to i$ and $\ell\to k$ and the four vertices are distinct,
then \eqref{frac:eq:edge-defect} gives
\((\eps_i+\eps_j)+(\eps_k+\eps_\ell)\geq2,\)
again contradicting $e<2$.

Finally, the arrows are transitive on three pairwise distinct vertices.
Suppose that
\(j\to i\) and \(k\to j\),
where $i,j,k$ are pairwise distinct. Set
\[
 \alpha:=m-\ee_j+\ee_i,
 \qquad
 \beta:=m-\ee_k+\ee_j.
\]
By the definition of the arrows, $\alpha,\beta\in\calS$. Their
coordinates satisfy
\[
 \alpha_i=\beta_i+1,
 \qquad
 \alpha_j=\beta_j-2,
 \qquad
 \alpha_k=\beta_k+1,
\]
and
\[
 \alpha_\ell=\beta_\ell
 \qquad
 \text{for every }\ell\notin\{i,j,k\}.
\]
In particular $\alpha_i>\beta_i$, and $j$ is the unique index $\ell$
for which $\alpha_\ell<\beta_\ell$. Apply the symmetric exchange axiom
to the ordered pair $(\alpha,\beta)$ at the coordinate $i$. Since the
exchange coordinate must satisfy $\alpha_\ell<\beta_\ell$, it is forced
to be $j$. Hence
\(\alpha-\ee_i+\ee_j=m\in\calS\)
and
\(\beta+\ee_i-\ee_j
 =m-\ee_k+\ee_i\in\calS.\)
Since $m_k>0$ by the arrow $k\to j$, the last membership says precisely
that $k\to i$. Thus the arrows of $G_m$ are transitive on every triple
of pairwise distinct vertices.

\smallskip\noindent
\emph{Step 2: the only possible nonempty graph configurations.}
Assume for contradiction that $G_m$ has an arrow, and fix one,
\(q\to p.\)
Let $r\in[s]\setminus\{p,q\}$ be a third vertex. There is no arrow
between $q$ and $r$. Indeed, an arrow $q\to r$ would give two distinct
outgoing arrows from $q$. On the other hand, if $r\to q$, then
transitivity applied to
\(r\to q\to p\)
would give $r\to p$, so $r$ would have the two distinct outgoing arrows
$r\to q$ and $r\to p$. Both possibilities are impossible.

Since two arrows cannot have disjoint endpoint sets, every other arrow
involving a third vertex must therefore meet $p$. Such an arrow cannot
be oriented as $p\to r$, because transitivity applied to
$q\to p\to r$ would give the forbidden arrow $q\to r$. Hence every
arrow involving a third vertex is of the form
\(r\to p.\)
There remains the possibility of the reciprocal arrow $p\to q$. If it
occurs and if there is also an arrow $r\to p$ from a third vertex, then
transitivity gives $r\to q$, so $r$ again has two distinct outgoing
arrows. Consequently there are only the following two possible
nontrivial configurations:
\begin{enumerate}[label=(\alph*)]
\item the only nonisolated vertices are $p$ and $q$, and
      $p\leftrightarrow q$;
\item every arrow points to a single centre vertex $p$, and $p$ has no
      outgoing arrow.
\end{enumerate}
The second case includes the possibility that $q\to p$ is the only
arrow.

We shall also need the following consequence for isolated vertices.
If $i\in[s]$ is isolated, meaning that no arrow is incident to $i$, then
\begin{equation}\label{frac:eq:isolated-fixed}
 m_i>0,
 \qquad
 \eps_i=0,
 \qquad
 \alpha_i=m_i\quad\text{for every }\alpha\in\calS.
\end{equation}
To prove this, first suppose that $m_i=0$. Put
\(A:=\sum_{r=1}^sL_r.\)
The divisor $A$ is ample, while $L_i$ is a nonzero nef Cartier divisor,
so Lemma~\ref{lem:ample-sum} gives
\(L_iA^{n-1}>0.\)
After expanding $A^{n-1}$, all mixed intersections are nonnegative, so
there is some $\beta\in\calS$ with
\(\beta_i>0=m_i.\)
Apply symmetric exchange to the ordered pair $(\beta,m)$ at the
coordinate $i$. There is an index $j$ with $\beta_j<m_j$ such that
\(m-\ee_j+\ee_i\in\calS.\)
Here $j\ne i$ and $m_j>0$, so this gives an incoming arrow $j\to i$,
contrary to isolation. Hence $m_i>0$.

Since $i$ has no outgoing arrow, the definition of the graph gives
\(I_{ki}=0\) for every \(k\ne i\).
Therefore \eqref{frac:eq:full-epsilon} gives immediately
\[
 \eps_i
 =\sum_{k\ne i}(1-\eps_k)I_{ki}
 =0.
\]

Finally let $\alpha\in\calS$. If $\alpha_i>m_i$, symmetric exchange
applied to $(\alpha,m)$ at $i$ produces an index $j$ with
\(m-\ee_j+\ee_i\in\calS,\)
and hence an incoming arrow $j\to i$. If instead $\alpha_i<m_i$,
symmetric exchange applied to $(m,\alpha)$ at $i$ produces an index
$j$ with
\(m-\ee_i+\ee_j\in\calS,\)
and hence an outgoing arrow $i\to j$. Both are impossible. Therefore
$\alpha_i=m_i$ for every $\alpha\in\calS$, proving
\eqref{frac:eq:isolated-fixed}.

\smallskip\noindent
\emph{Step 3: exclusion of the reciprocal pair.}
Consider case (a). Thus
\(q\to p\) and \(p\to q\),
and there are no other arrows. In particular $m_p,m_q>0$. Put
\(h:=I_{pq}\) and \(k:=I_{qp}\).
Since all vertices other than $p,q$ are isolated, the only nonzero
off-diagonal term in \eqref{frac:eq:full-epsilon} for the source $q$
is $I_{pq}=h$, and the only such term for the source $p$ is
$I_{qp}=k$. Hence
\begin{equation}\label{frac:eq:two-cycle-epsilon}
 \eps_q=h(1-\eps_p),
 \qquad
 \eps_p=k(1-\eps_q).
\end{equation}

The mixed Khovanskii--Teissier inequality forces $h=k=1$. Indeed, set
\[
 N:=L_p^{m_p-1}L_q^{m_q-1}
 \prod_{i\notin\{p,q\}}L_i^{m_i}.
\]
Then
\[
 L_pL_qN=L^m=1,
 \qquad
 L_p^2N=I_{pq}=h,
 \qquad
 L_q^2N=I_{qp}=k.
\]
Applying \eqref{eq:mixed-KT} to $L_p,L_q$ and the $(n-2)$-fold product
$N$ gives
\(hk=(L_p^2N)(L_q^2N)
 \leq(L_pL_qN)^2
 =1.\)
Since $h$ and $k$ are positive integers, $h=k=1$. Substituting in
\eqref{frac:eq:two-cycle-epsilon} gives
\begin{equation}\label{frac:eq:two-cycle-defect}
 \eps_p+\eps_q=1.
\end{equation}
In particular $0<\eps_p,\eps_q<1$.

Put
\(\mu:=m_p+m_q\)
and define
\[
 b_\ell
 :=L_p^\ell L_q^{\mu-\ell}
   \prod_{i\notin\{p,q\}}L_i^{m_i}
 \qquad(0\leq\ell\leq\mu).
\]
By Lemma~\ref{lem:exchange}, the positive support of the sequence
$(b_\ell)$ is an interval
\([u,v]\cap\ZZ.\)
Every vertex outside $\{p,q\}$ is isolated, so
\eqref{frac:eq:isolated-fixed} gives
\[
 \alpha_i=m_i
 \qquad
 \text{for every }\alpha\in\calS
 \text{ and every }i\notin\{p,q\}.
\]
Since every $\alpha\in\calS$ has total degree $n=|m|$, it follows that
\(\alpha_p+\alpha_q=m_p+m_q=\mu.\)
By \eqref{frac:eq:isolated-fixed}, every coordinate outside
$\{p,q\}$ is fixed equal to $m_i$ on $\calS$. Since
$|\alpha|=|m|=n$, every $\alpha\in\calS$ therefore satisfies
$\alpha_q=\mu-\alpha_p$. Hence
\(L^\alpha
 =
 b_{\alpha_p},\)
and conversely every $\ell$ with $b_\ell>0$ determines an element of
$\calS$ by taking
$\alpha_p=\ell$, $\alpha_q=\mu-\ell$, and
$\alpha_i=m_i$ for $i\notin\{p,q\}$.
Consequently
\begin{equation}\label{frac:eq:two-cycle-support}
 \left\{\alpha_p\ \middle|\ \alpha\in\calS\right\}
 =\left\{\ell\in\{0,\ldots,\mu\}\ \middle|\  b_\ell>0\right\}
 =[u,v]\cap\ZZ.
\end{equation}
In particular,
\[
 u=\min_{\alpha\in\calS}\alpha_p,
 \qquad
 v=\max_{\alpha\in\calS}\alpha_p.
\]
Moreover,
\[
 b_{m_p-1}=I_{qp}=1,
 \qquad
 b_{m_p}=L^m=1,
 \qquad
 b_{m_p+1}=I_{pq}=1.
\]

We claim that $u=0$. Suppose that $u>0$, and set
$B_p:=\sum_{i\ne p}L_i$.
Using \eqref{eq:numdim} and \eqref{frac:eq:two-cycle-support}, we obtain
\(\nu(B_p)=\max_{\alpha\in\calS}\sum_{i\ne p}\alpha_i =\max_{\alpha\in\calS}(n-\alpha_p)=n-u.\)
The nef Cartier divisor $B_p$ is semiample. Let
$g_p:X\longrightarrow Y_p$ be its contraction. By Lemma~\ref{lem:numdim-basic},
$\dim Y_p=n-u$, so a general fibre $F$ has dimension $u$.
Since $u>0$, this is a positive-dimensional general fibre and hence is a
klt Fano variety. On $F$ one has
\(B_p|_F\equiv0.\)
Every $L_i$ with $i\ne p$ is nef and is a summand of $B_p$, so for every
curve $C\subset F$,
\[
 0=B_p\cdot C=\sum_{i\ne p}L_i\cdot C
\]
with all terms nonnegative. Thus
\(L_i|_F\equiv0\) for \(i\ne p\).
Restricting the anticanonical decomposition to $F$ therefore gives
$-K_F\equiv(m_p+1-\eps_p)L_p|_F$.
Because $0<\eps_p<1$, the rounding integer associated with this
one-term decomposition is exactly $m_p$. Applying
Lemma~\ref{lem:Bm-vanishing} on the $u$-dimensional klt Fano variety
$F$ gives
\(m_p\leq u.\)
On the other hand $b_{m_p-1}=1$, so
$m_p-1\in[u,v]$ and hence
\(u\leq m_p-1.\)
This is impossible. Therefore $u=0$.

The same argument with $p$ and $q$ interchanged gives $v=\mu$. More
explicitly, the minimum value of $\alpha_q$ on $\calS$ is
$\mu-v$. If $v<\mu$, applying the preceding fibre argument to
$B_q:=\sum_{i\ne q}L_i$ gives
\(m_q\leq\mu-v.\)
But $b_{m_p+1}=1$ corresponds to a support vector with
$q$-coordinate $m_q-1$, so
\(\mu-v\leq m_q-1,\)
a contradiction. Hence $v=\mu$.

Thus $b_\ell>0$ for every $0\leq\ell\leq\mu$. The sequence consists of
positive integers and is log-concave:
\(b_\ell^2\geq b_{\ell-1}b_{\ell+1}.\)
Starting from
\(b_{m_p-1}=b_{m_p}=b_{m_p+1}=1,\)
log-concavity propagates the value one in both directions. For example,
\[
 1=b_{m_p-1}^2
 \geq b_{m_p-2}b_{m_p}=b_{m_p-2},
\]
so positivity and integrality give $b_{m_p-2}=1$. Repeating this argument
moves successively to the left. The same argument starting from
$b_{m_p+1}=b_{m_p}=1$ moves to the right. Therefore
\begin{equation}\label{frac:eq:all-b-one}
 b_0=b_1=\cdots=b_\mu=1.
\end{equation}

We now compute the volume polynomial explicitly. By
\eqref{frac:eq:isolated-fixed}, a multi-index $\alpha\in\calS$ has the
form
\[
 \alpha_p=\ell,
 \qquad
 \alpha_q=\mu-\ell,
 \qquad
 \alpha_i=m_i\quad(i\notin\{p,q\}),
\]
for some $0\leq\ell\leq\mu$, and
\eqref{frac:eq:all-b-one} says that the corresponding mixed
intersection equals one. Therefore the multinomial expansion
\eqref{eq:volume-polynomial} gives
\[
 \begin{aligned}
 \vol_L(t)
 &=\sum_{\ell=0}^{\mu}
 \frac{n!}{\ell!(\mu-\ell)!\prod_{i\notin\{p,q\}}m_i!}
 t_p^\ell t_q^{\mu-\ell}
 \prod_{i\notin\{p,q\}}t_i^{m_i}\\
 &=\frac{n!}{\mu!\prod_{i\notin\{p,q\}}m_i!}
 \left(\sum_{\ell=0}^{\mu}\binom{\mu}{\ell}
 t_p^\ell t_q^{\mu-\ell}\right)
 \prod_{i\notin\{p,q\}}t_i^{m_i}\\
 &=\frac{n!}{\mu!\prod_{i\notin\{p,q\}}m_i!}
 (t_p+t_q)^\mu
 \prod_{i\notin\{p,q\}}t_i^{m_i}.
 \end{aligned}
\]
Set $A:=\sum_{i=1}^sL_i$ and $D:=L_p-L_q$.
We now use the volume polynomial to show that the two classes
$L_p$ and $L_q$ are numerically indistinguishable.  This will
contradict the linear independence of the optimal vertex decomposition.
Evaluate the volume polynomial along the one-parameter line
$t_p=1+z$, $t_q=1-z$, and $t_i=1$ for $i\notin\{p,q\}$.
Since $t_p+t_q=2$, the displayed formula for $\vol_L(t)$ is independent
of $z$. On the other hand, by definition of the volume polynomial this
specialization is
\((A+zD)^n.\)
Hence $(A+zD)^n$ is constant as a polynomial in $z$. Its first two
derivatives at $z=0$ give
\(DA^{n-1}=0\) and \(D^2A^{n-2}=0\).
The divisor $A$ is ample, so the numerical Hodge index theorem
\eqref{eq:hodge} yields $D\equiv0$. Thus
\(L_p\equiv L_q,\)
contradicting the numerical linear independence of the classes
$[L_1],\ldots,[L_s]$ in an optimal vertex decomposition. This excludes
case (a).

\smallskip\noindent
\emph{Step 4: exclusion of the inward star.}
Consider case (b). Thus every arrow points to a single centre
vertex $p$, and $p$ has no outgoing arrow. Since the graph is nonempty,
$p$ is an endpoint of an arrow. By \eqref{frac:eq:edge-defect},
\(\eps_p>0.\)
We first show that
\(m_p=0.\)
If $m_p>0$, equation \eqref{frac:eq:full-epsilon} applies to the source
vertex $p$. Since $\eps_p>0$, at least one term on the right-hand side
would be positive, so $I_{ip}>0$ for some $i\ne p$. This would give an
outgoing arrow $p\to i$, contrary to the definition of case
(b). Hence $m_p=0$.

Every vertex $i\ne p$ satisfies $m_i>0$. If $i$ is nonisolated, then it
is the source of an arrow $i\to p$, and the definition of an arrow gives
$m_i>0$. If $i$ is isolated, this follows from
\eqref{frac:eq:isolated-fixed}.

Put
\[
 b:=1-\eps_p\in(0,1),
 \qquad
 h_i:=I_{pi}\in\ZZ_{\geq0}\quad(i\ne p),
 \qquad
 h:=\sum_{i\ne p}h_i.
\]
For $i\ne p$, the only possible outgoing arrow is $i\to p$. Thus in
\eqref{frac:eq:full-epsilon} with source $i$, every off-diagonal term
vanishes except possibly the term indexed by $p$. Consequently
\begin{equation}\label{frac:eq:star-epsilon}
 \eps_i=(1-\eps_p)I_{pi}=bh_i
 \qquad(i\ne p).
\end{equation}
Since the graph is nonempty, at least one $h_i$ is positive, so $h\geq1$.
Moreover $bh_i=\eps_i<1$ for every $i\ne p$. Summing
\eqref{frac:eq:star-epsilon} gives
\[
 e=\eps_p+\sum_{i\ne p}\eps_i
   =(1-b)+bh
   =1+b(h-1)<2.
\]

We next determine the numerical dimensions of the classes $L_i$ for
$i\ne p$. No arrow enters such a vertex $i$. If some
$\alpha\in\calS$ satisfied $\alpha_i>m_i$, apply symmetric exchange to
$(\alpha,m)$ at the coordinate $i$. It produces an index $j$ with
$\alpha_j<m_j$ and
\(m-\ee_j+\ee_i\in\calS.\)
Then $m_j>0$ and this is an incoming arrow $j\to i$, a contradiction.
Therefore
\[
 \alpha_i\leq m_i
 \qquad
 \text{for every }\alpha\in\calS.
\]
Since $m\in\calS$, equality is attained, and formula \eqref{eq:numdim}
with $J=\{i\}$ gives
\begin{equation}\label{frac:eq:star-numdim}
 \nu(L_i)=\max_{\alpha\in\calS}\alpha_i=m_i
 \qquad(i\ne p).
\end{equation}

Let
\(f_i:X\longrightarrow Y_i\)
be the contraction defined by the semiample divisor $L_i$. By
\eqref{frac:eq:star-numdim}, $\dim Y_i=m_i$. For every integer
$1\leq t\leq m_i$, the vector $t\ee_i$ lies in the box
$\mathcal B_m$, and Lemma~\ref{lem:Bm-vanishing} gives
\[
 \chi(X,-tL_i)=0.
\]
Lemma~\ref{lem:image-P} therefore identifies the polarized image with
projective space:
\[
 (Y_i,H_i)\simeq(\PP^{m_i},\OO(1)),
 \qquad
 L_i\sim f_i^*H_i.
\]
After this identification, write simply
\[
 f_i:X\longrightarrow\PP^{m_i},
 \qquad
 L_i=f_i^*\OO_{\PP^{m_i}}(1).
\]
Consider the product morphism
\[
 f:=(f_i)_{i\ne p}:
 X\longrightarrow
 Y:=\prod_{i\ne p}\PP^{m_i}.
\]
Since $m_p=0$ and $M=n$, $\dim Y=\sum_{i\ne p}m_i=n=\dim X$.
Furthermore, $\prod_{i\ne p}L_i^{m_i}=L^m=1$.
If the image of $f$ had dimension less than $n$, this top-degree
intersection of divisor classes pulled back from $Y$ would vanish.
Hence $f$ is dominant, and therefore generically finite. By the projection
formula,
\[
 1=\prod_{i\ne p}L_i^{m_i}
  =\deg(f)\prod_{i\ne p}H_i^{m_i}
  =\deg(f).
\]
Thus $f$ is birational.

Finally, using $m_p=0$, $b=1-\eps_p$, and
\eqref{frac:eq:star-epsilon}, the anticanonical decomposition becomes
\[
 -K_X
 \equiv
 bL_p+\sum_{i\ne p}(m_i+1-bh_i)L_i.
\]
We have
\[
 0<b<1,
 \qquad
 h_i\in\ZZ_{\geq0},
 \qquad
 bh_i<1,
 \qquad
 \sum_{i\ne p}h_i=h\geq1,
\]
and the birational morphism $f$ has target
$\prod_{i\ne p}\PP^{m_i}$ with $\sum_{i\ne p}m_i=n$.
Thus all hypotheses of Lemma~\ref{frac:lem:star} are satisfied, with
$L=L_p$. That lemma gives
\[
 L_p\equiv\sum_{i\ne p}h_iL_i.
\]
This contradicts the numerical linear independence of
$[L_1],\ldots,[L_s]$. Hence case (b) is also impossible.

\smallskip\noindent
\emph{Step 5: conclusion.}
Both possible nonempty graph configurations have been excluded. Therefore
$G_m$ has no arrows. Equivalently,
\[
 I_{ij}=0
 \qquad
 \text{for every $j$ with $m_j>0$ and every $i\ne j$}.
\]
Lemma~\ref{lem:singleton} now gives
\[
 m_i>0\quad\text{for every }i \quad\text{and}\quad
 \calS=\{m\}.
\]
Since every $m_j>0$, equation \eqref{frac:eq:full-epsilon} applies to
every $j$ and its right-hand side is zero. Hence
\(\eps_j=0\) for every \(j\),
so $e=0$. The anticanonical decomposition is therefore the integral
full-degree decomposition
\[
 -K_X\equiv\sum_{i=1}^s(m_i+1)L_i,
 \qquad
 \sum_i m_i=n.
\]
Proposition~\ref{prop:product} yields
\[
 X\simeq\prod_{i=1}^s\PP^{m_i}
 \qquad\text{and}\qquad
 \rho(X)=s.
\]
Finally,
\[
 c_X=(n-M)+(\rho(X)-s)+e=0.
\]
This proves the proposition.
\end{proof}

\begin{corollary}\label{frac:cor:reduction}
If $1\leq c_X<2$, then every optimal vertex decomposition satisfies
\begin{equation}\label{frac:eq:reduction}
 M=n-1,\qquad s=\rho(X),\qquad0\leq e=c_X-1<1.
\end{equation}
\end{corollary}
\begin{proof}
The nonnegative integral terms $n-M$ and $\rho-s$ in
\eqref{eq:defect} have sum at most one. The possibility $M=n$ is excluded
by Proposition~\ref{frac:prop:full}. The remaining possibility is
\eqref{frac:eq:reduction}.
\end{proof}

\section{Inactive factors and product decompositions}\label{frac:sec:mixed}

Throughout this section, let $X$ be a complex klt Fano variety of
dimension $n$ with $1\leq c_X<2$, and fix one optimal vertex decomposition. By
Corollary~\ref{frac:cor:reduction},
\begin{equation}\label{frac:eq:mixed-data}
 -K_X\equiv\sum_{i=1}^{s}(m_i+1-\eps_i)L_i,
 \qquad |m|=n-1,\quad s=\rho(X),\quad0\leq e=c_X-1<1.
\end{equation}
Here and below $m=(m_1,\ldots,m_s)$.  We retain from
Section~\ref{sec:equality} the notation
$L^\alpha=\prod_iL_i^{\alpha_i}$,
$\calS=\{\alpha\in\ZZ_{\geq0}^s\mid|\alpha|=n,\ L^\alpha>0\}$, the associated
Euler polynomial $P(k)=\chi(X,\OO_X(-\sum_i k_iL_i))$, and
$\Delta_iP(t)=P(t+\ee_i)-P(t)$,
$\Delta^u=\Delta_1^{u_1}\cdots\Delta_s^{u_s}$.  Set
\[
 d_i=L^{m+\ee_i},\qquad I=\left\{i\ \middle|\  d_i>0\right\},\qquad
 J=\left\{i\ \middle|\  d_i=0\right\}.
\]
We call the indices in $I$ \emph{active} and those in $J$ \emph{inactive}.
We record explicitly the numerical identity that drives the remainder of
the argument.  Since $|m|=n-1$, the vanishing-box property of
Lemma~\ref{lem:Bm-vanishing} and the finite-difference formula give
\[
 \Delta^mP(0)=(-1)^{n-1}.
\]
Applying \eqref{eq:RR-difference} with $u=m$ therefore yields
\[
 \left(K_X+\sum_i m_iL_i\right)L^m
 =2(-1)^n\Delta^mP(0)=-2.
\]
On the other hand, \eqref{frac:eq:mixed-data} gives
\[
 K_X+\sum_i m_iL_i
 \equiv-\sum_i(1-\eps_i)L_i.
\]
Intersecting with $L^m$ and using
$L_iL^m=L^{m+\ee_i}=d_i$, we obtain
\begin{equation}\label{frac:eq:weighted}
 \sum_i(1-\eps_i)d_i=2.
\end{equation}
Every coefficient $1-\eps_i$ is positive and every $d_i$ is a
nonnegative integer.  Hence at least one $d_i$ is positive, so
$I\ne\varnothing$.

\subsection{The inactive base}

\begin{lemma}\label{frac:lem:tightness}
One has
\[
 \nu\left(\sum_{j\in J}L_j\right)=\sum_{j\in J}m_j.
\]
\end{lemma}
\begin{proof}
For a subset $T\subseteq[s]$, write $m(T):=\sum_{i\in T}m_i$.
We call $T$ \emph{tight} if
\(r(T)=m(T).\)
Let $r$ be the polymatroid rank function of \eqref{eq:rank}, so that
\[
 r(T)=\max_{\alpha\in\calS}\alpha(T)
     =\nu\left(\sum_{i\in T}L_i\right).
\]
We first establish the basic lower bound
\begin{equation}\label{frac:eq:m-below-r}
 m(T)\leq r(T)
 \qquad\text{for every }T.
\end{equation}
Indeed, choose any active index $p\in I$.  Then
$d_p=L^{m+\ee_p}>0$, so $m+\ee_p\in\calS$.  Hence
\(r(T)\geq(m+\ee_p)(T)\geq m(T).\)

Now fix an inactive index $j\in J$.  Since $d_j=0$, the vector
$m+\ee_j$ does not belong to $\calS$.  In fact there is no
$\alpha\in\calS$ with $m+\ee_j\leq\alpha$ coordinatewise: both vectors
have total degree $n$, so such an inequality would force
$\alpha=m+\ee_j$.  By the rank-extension criterion
\eqref{eq:rank-extension}, there is therefore a subset $T_j$ such that
\begin{equation}\label{frac:eq:Tj-violation}
 (m+\ee_j)(T_j)>r(T_j).
\end{equation}
The set $T_j$ must contain $j$.  Otherwise the left-hand side of
\eqref{frac:eq:Tj-violation} would be $m(T_j)$, contradicting
\eqref{frac:eq:m-below-r}.  Since $j\in T_j$,
\((m+\ee_j)(T_j)=m(T_j)+1.\)
Together with \eqref{frac:eq:m-below-r}, this gives
\(m(T_j)\leq r(T_j)<m(T_j)+1.\)
The rank $r(T_j)$ is an integer, hence
\begin{equation}\label{frac:eq:Tj-tight}
 r(T_j)=m(T_j).
\end{equation}
Moreover $T_j$ contains no active index.  If $p\in I\cap T_j$, then
$m+\ee_p\in\calS$ would give
\(r(T_j)\geq(m+\ee_p)(T_j)=m(T_j)+1,\)
contrary to \eqref{frac:eq:Tj-tight}.  Thus
\begin{equation}\label{frac:eq:Tj-in-J}
 j\in T_j\subseteq J.
\end{equation}

We next note that unions of tight sets are tight.  Suppose
$r(T)=m(T)$ and $r(U)=m(U)$.  Submodularity gives
\[
 r(T\cup U)+r(T\cap U)\leq r(T)+r(U)=m(T)+m(U).
\]
Since $m$ is additive on unions and intersections,
\(m(T)+m(U)=m(T\cup U)+m(T\cap U).\)
The lower bound \eqref{frac:eq:m-below-r} applies separately to
$T\cup U$ and $T\cap U$.  Therefore equality must hold in both terms,
and in particular
\(r(T\cup U)=m(T\cup U).\)
By \eqref{frac:eq:Tj-in-J}, the union of the sets $T_j$, $j\in J$, is
exactly $J$.  Repeatedly taking unions therefore yields
\[
 r(J)=m(J)=\sum_{j\in J}m_j.
\]
Finally, the numerical-dimension formula \eqref{eq:rank} gives
\[
 \nu\left(\sum_{j\in J}L_j\right)=r(J)=\sum_{j\in J}m_j,
\]
as required.
\end{proof}

\begin{lemma}\label{frac:lem:inactive-base}
If $J\ne\varnothing$, the contraction of $\sum_{j\in J}L_j$ is
\[
 f:X\longrightarrow Y=\prod_{j\in J}\PP^{m_j},
 \qquad L_j\sim f^*H_j\quad(j\in J),
\]
where all $m_j>0$. In particular, $\nu(L_j)=m_j$ for every inactive $j$.
For $J=\varnothing$ take $Y$ to be a point.
\end{lemma}
\begin{proof}
The case $J=\varnothing$ is part of the convention in the statement, so
assume $J\ne\varnothing$ and put $D_J:=\sum_{j\in J}L_j$.
The divisor $D_J$ is nef Cartier.  Since $X$ is klt Fano,
Lemma~\ref{lem:nef-monoid} shows that $D_J$ is semiample.  Let
\(f:X\longrightarrow Y\)
be the contraction defined by $D_J$, with $Y$ normal and the fibres of
$f$ connected.

\smallskip\noindent
\emph{Step 1: descent of the inactive divisors.}
Let $C$ be a curve contracted by $f$.  Then
\[
 0=D_JC=\sum_{j\in J}L_jC.
\]
Every summand is nonnegative because the $L_j$ are nef.  Hence
$L_jC=0$ for every $j\in J$.  Thus each $L_j$ is numerically trivial
over $Y$.  Since $-K_X$ is $f$-ample, Lemma~\ref{lem:descent} gives a
Cartier divisor $H_j$ on $Y$ such that
\begin{equation}\label{frac:eq:inactive-descent}
 L_j\sim f^*H_j.
\end{equation}
Each $H_j$ is nef: if $B\subset Y$ is an integral curve, choose an
integral curve $C\subset X$ mapping generically finitely onto $B$; then
\(L_jC=\deg(C/B)\,H_jB\geq0.\)
No $H_j$ is numerically trivial, since otherwise
$L_j=f^*H_j\equiv0$, contrary to the choice of the supporting classes.
Moreover
\[
 D_J\sim f^*\left(\sum_{j\in J}H_j\right).
\]
Because $f$ is the contraction defined by $D_J$, the final assertion of
Lemma~\ref{lem:descent} shows that
\begin{equation}\label{frac:eq:inactive-ample-sum}
 \sum_{j\in J}H_j
 \quad\text{is ample on }Y.
\end{equation}

Since $D_J$ is semiample,
Lemma~\ref{lem:numdim-basic} identifies its numerical dimension with the
dimension of its image.  Lemma~\ref{frac:lem:tightness} therefore gives
\begin{equation}\label{frac:eq:inactive-base-dimension}
 \dim Y=\nu(D_J)=\sum_{j\in J}m_j.
\end{equation}

\smallskip\noindent
\emph{Step 2: move the active part into a boundary.}
For brevity put $a_i:=m_i+1-\eps_i$.
Every active divisor $L_i$, $i\in I$, is nef Cartier and hence
semiample.  Choose a common sufficiently divisible integer $N$ such that
all $|NL_i|$, $i\in I$, are base-point-free and such that
$a_i/N<1$ for every active $i$.  Choose general members
\(D_i\in|NL_i|\) for \(i\in I\)
and set $\Delta:=\sum_{i\in I}\frac{a_i}{N}D_i$.
Since $X$ is klt, and since the systems $|NL_i|$ are base-point-free
with $a_i/N<1$ for every active $i$, Lemma~4.7(1) and
Theorem~4.8(2) of \cite{KollarPairs} imply that
\begin{equation}\label{frac:eq:Delta-klt}
 (X,\Delta)\quad\text{is klt}.
\end{equation}
Furthermore $\Delta\sim_{\QQ}\sum_{i\in I}a_iL_i$.
On a klt Fano variety numerical equivalence of $\QQ$-Cartier divisors
agrees with $\QQ$-linear equivalence.  Hence the optimal decomposition
\eqref{frac:eq:mixed-data} and \eqref{frac:eq:inactive-descent} give
\begin{equation}\label{frac:eq:inactive-lc-trivial}
 K_X+\Delta\sim_{\QQ}-\sum_{j\in J}a_jL_j
 \sim_{\QQ}f^*\left(-\sum_{j\in J}a_jH_j\right).
\end{equation}
Thus $K_X+\Delta$ is $\QQ$-linearly trivial over $Y$.

\smallskip\noindent
\emph{Step 3: the canonical bundle formula on the base.}
Apply the canonical bundle formula
\cite[Theorem~3.1]{FG} to the klt pair $(X,\Delta)$ and the morphism
$f$.  It yields an effective $\QQ$-divisor $\Gamma$ on $Y$ such that
$(Y,\Gamma)$ is klt and
\[
 K_X+\Delta\sim_{\QQ}f^*(K_Y+\Gamma).
\]
Comparing this with \eqref{frac:eq:inactive-lc-trivial}, we obtain
\[
 f^*\left(K_Y+\Gamma+\sum_{j\in J}a_jH_j\right)
 \sim_{\QQ}0.
\]
The pullback map on $\QQ$-Picard groups is injective.  Indeed,
$f_*\OO_X=\OO_Y$ because the fibres are connected and $Y$ is normal;
if a Cartier divisor $B$ satisfies $f^*B\sim0$, then the projection
formula gives
\[
 \OO_Y(B)
 \simeq f_*f^*\OO_Y(B)
 \simeq f_*\OO_X
 \simeq\OO_Y.
\]
After clearing denominators, the same argument applies to
$\QQ$-Cartier divisors.  Consequently
\begin{equation}\label{frac:eq:inactive-log-Fano-base}
 -(K_Y+\Gamma)
 \sim_{\QQ}
 \sum_{j\in J}(m_j+1-\eps_j)H_j.
\end{equation}
The right-hand side is ample.  Indeed, all coefficients are positive,
all $H_j$ are nef, and their sum is ample by
\eqref{frac:eq:inactive-ample-sum}.

\smallskip\noindent
\emph{Step 4: the full-degree classification of the base.}
Equation~\eqref{frac:eq:inactive-base-dimension} says exactly that
\[
 \sum_{j\in J}m_j=\dim Y.
\]
Moreover
\[
 \sum_{j\in J}\eps_j\leq\sum_i\eps_i=e<1.
\]
Thus all hypotheses of Proposition~\ref{prop:product} apply to the klt
log Fano pair $(Y,\Gamma)$ and the decomposition
\eqref{frac:eq:inactive-log-Fano-base}.  We conclude that every
$m_j$ is positive and that
\[
 (Y;H_j)_{j\in J}
 \simeq
 \left(
 \prod_{j\in J}\PP^{m_j};
 p_j^*\OO_{\PP^{m_j}}(1)
 \right)_{j\in J}.
\]
Together with \eqref{frac:eq:inactive-descent}, this gives
\(L_j\sim f^*p_j^*\OO_{\PP^{m_j}}(1).\)
Finally $L_j$ is semiample and its associated image is the $j$th
projective-space factor, so Lemma~\ref{lem:numdim-basic} gives
\[
 \nu(L_j)=m_j.
\]
This proves all assertions.
\end{proof}

\subsection{Integral inactive factors}

\begin{lemma}\label{frac:lem:inactive-split}
Let $J_0=\left\{j\in J\ \middle|\ \eps_j=0\right\}$. Then
\[
 X\simeq Z\times\prod_{j\in J_0}\PP^{m_j},
\]
where $Z$ is klt Fano, $c_Z=c_X$, and
$\rho(Z)=\rho(X)-|J_0|$. The classes $L_j$, $j\in J_0$, are the
hyperplane classes of the displayed factors, and all other $L_i$ descend
to $Z$. The descended decomposition is optimal and independent.
\end{lemma}
\begin{proof}
If $J_0=\varnothing$, take $Z=X$ and there is nothing to prove.
Assume henceforth that $J_0\ne\varnothing$. We first record a support
statement for every $j\in J_0$.

Fix $j\in J_0$ and put
\(B_j:=\sum_{i\ne j}L_i.\)
We claim that $B_j$ is not big.  Suppose otherwise.  The vanishing-box property
for the optimal decomposition is valid on
\[
 \mathcal B_m=\left\{k\ \middle|\  0\leq k_i\leq m_i\right\}.
\]
Because $\eps_j=0$, consider the enlarged box
\[
 \mathcal B_m^{(j)}
 :=
 \left\{
 k\in\ZZ_{\geq0}^s
 \ \middle|\
 0\leq k_i\leq m_i\ (i\ne j),\quad
 0\leq k_j\leq m_j+1
 \right\}.
\]
The vanishing is already known on the subbox $\mathcal B_m$.
It therefore remains only to consider the new face
\[
 k_j=m_j+1,
 \qquad
 0\leq k_i\leq m_i\quad(i\ne j).
\]
For such $k$, since $\eps_j=0$,
\[
 -K_X-\sum_i k_iL_i
 \equiv
 \sum_{i\ne j}(m_i+1-\eps_i-k_i)L_i.
\]
For every $i\ne j$,
\(m_i+1-\eps_i-k_i\geq1-\eps_i>0.\)
Hence, under the assumption that
$B_j=\sum_{i\ne j}L_i$ is big, the divisor on the right is nef and big.
Kawamata--Viehweg vanishing and the same non-effectivity argument as in
Lemma~\ref{lem:Bm-vanishing} therefore give
\[
 P(0)=1,
 \qquad
 P(k)=0
 \quad\text{for every nonzero }k\in\mathcal B_m^{(j)}.
\]
Put $u=m+\ee_j$.  Since $|u|=n$, the finite-difference formula gives
\[
 \Delta^uP(0)=(-1)^n.
\]
On the other hand Lemma~\ref{lem:finite-differences} (iii), applied to the top Riemann--Roch term, gives
\[
 \Delta^uP(0)=(-1)^nL^u=(-1)^nd_j=0,
\]
a contradiction.  Hence
\begin{equation}\label{frac:eq:Bj-not-big}
 B_j\quad\text{is not big.}
\end{equation}

Let $h_j:X\to W_j$ be the contraction of $B_j$, and let $F_j$ be a
general fibre.  Put
\(r_j:=\dim F_j.\)
By \eqref{frac:eq:Bj-not-big}, $r_j>0$.  The numerical-dimension formula
\eqref{eq:numdim} gives
\[
 r_j=n-\nu(B_j)
     =\min_{\alpha\in\calS}\alpha_j.
\]
Every $L_i$ with $i\ne j$ is numerically trivial on $F_j$, while
$L_j|_{F_j}$ is ample because $\sum_iL_i$ is ample.  Since
$\eps_j=0$, adjunction on a general fibre gives
\[
 -K_{F_j}\equiv(m_j+1)L_j|_{F_j}.
\]
Lemma~\ref{lem:index} therefore gives
\begin{equation}\label{frac:eq:rj-lower}
 r_j\geq m_j.
\end{equation}
Choose an active index $p\in I$.  Since
$m+\ee_p\in\calS$ and $p\ne j$,
\[
 r_j=\min_{\alpha\in\calS}\alpha_j
 \leq(m+\ee_p)_j=m_j.
\]
Together with \eqref{frac:eq:rj-lower}, this gives
$r_j=m_j$, and hence
\begin{equation}\label{frac:eq:support-lower-inactive}
 \alpha_j\geq m_j
 \qquad(\alpha\in\calS).
\end{equation}
On the other hand Lemma~\ref{frac:lem:inactive-base} gives
$\nu(L_j)=m_j$.  Applying \eqref{eq:numdim} to the singleton
$\{j\}$ gives
\(\max_{\alpha\in\calS}\alpha_j=m_j.\)
Combining this with \eqref{frac:eq:support-lower-inactive}, we obtain
\begin{equation}\label{frac:eq:fixed-inactive}
 \alpha_j=m_j
 \qquad(\alpha\in\calS,\ j\in J_0).
\end{equation}

Now put
\(D':=\sum_{i\notin J_0}L_i\)
and let
\(g:X\longrightarrow Z\)
be its contraction.  If a curve is contracted by $g$, every summand of
$D'$ has zero intersection with that curve, so Lemma~\ref{lem:descent}
shows that all $L_i$, $i\notin J_0$, descend to nef Cartier divisors on
$Z$.  Their sum is ample.  For every $\alpha\in\calS$, one has
$|\alpha|=n$, while \eqref{frac:eq:fixed-inactive} gives
$\alpha_j=m_j$ for every $j\in J_0$.  Hence
\(\sum_{i\notin J_0}\alpha_i =n-\sum_{j\in J_0}m_j,\)
independently of $\alpha\in\calS$.  Formula \eqref{eq:numdim} therefore gives
\[
 \dim Z=\nu(D')
 =\max_{\alpha\in\calS}\sum_{i\notin J_0}\alpha_i
 =n-\sum_{j\in J_0}m_j.
\]
For $j\in J_0$, Lemma~\ref{frac:lem:inactive-base} gives a morphism to
$\PP^{m_j}$ with $L_j$ equal to the pullback of the hyperplane class.
Combining all these morphisms gives
\[
 \Phi:X\longrightarrow
 Z\times\prod_{j\in J_0}\PP^{m_j}.
\]
The pullback of the sum of an ample divisor on $Z$ and the hyperplane
classes on all projective factors is $\sum_iL_i$, which is ample.
Consequently $\Phi$ has no positive-dimensional fibre.  It is projective,
so it is finite onto its image.  The target has dimension
$\dim Z+\sum_{j\in J_0}m_j=n$, therefore the finite image is the whole target.

We next compute the degree. A general fibre $F$ of $g$ is normal and
klt and has dimension $\dim F=\sum_{j\in J_0}m_j$.
The divisors $L_i$ with $i\notin J_0$ are numerically trivial on $F$,
and adjunction gives
\[
 -K_F\equiv
 \sum_{j\in J_0}(m_j+1)L_j|_F.
\]
The restriction of $\Phi$ to $F$ is finite and surjective onto
$\prod_{j\in J_0}\PP^{m_j}$. Thus the classes $L_j|_F$ are nonzero nef
Cartier divisors with ample sum, and the sum of the integers $m_j$
equals $\dim F$.  Proposition~\ref{prop:product}
therefore identifies the induced finite morphism
$F\longrightarrow\prod_{j\in J_0}\PP^{m_j}$ with a degree-one morphism.  Equivalently,
\[
 \prod_{j\in J_0}(L_j|_F)^{m_j}=1.
\]
We now compare the degree of this fibre map with the degree of $\Phi$.
Put $T:=Z\times\prod_{j\in J_0}\PP^{m_j}$.
Since $\Phi_*\OO_X$ is coherent and its generic stalk is the
$K(T)$-vector space $K(X)$ of dimension $\deg\Phi$,
\cite[Chapter~II, Exercise~5.7(a)]{Hartshorne} gives, after replacing
$T$ by a dense open subset $U\subseteq T$, that
\[
 \Phi^{-1}(U)\longrightarrow U
\]
is finite locally free of degree $\deg\Phi$.

Let $q:T\longrightarrow Z$ be the first projection. Since $U$ is a nonempty open subset of the
product $T$, its image $q(U)$ contains a nonempty open subset of $Z$.
Choose $z\in Z$ general in this open subset, so that the preceding
properties of the general fibre $F=g^{-1}(z)$ also hold. Then
\[
 U_z:=U\cap\left(\{z\}\times
 \prod_{j\in J_0}\PP^{m_j}\right)
\]
is a nonempty open subset of
$\{z\}\times\prod_{j\in J_0}\PP^{m_j}$. Since
$F=\Phi^{-1}(\{z\}\times\prod_{j\in J_0}\PP^{m_j})$, base change
preserves finite locally free morphisms and their rank. Hence
\[
 (\Phi|_F)^{-1}(U_z)\longrightarrow U_z
\]
is finite locally free of degree $\deg\Phi$, and therefore the finite morphism
\[
 \Phi|_F:F\longrightarrow\prod_{j\in J_0}\PP^{m_j}
\]
has degree $\deg\Phi$. The preceding intersection computation shows
that this fibre map has degree one. Therefore
\(\deg\Phi=1.\)
The target is normal.  Thus the finite birational morphism $\Phi$ is an
isomorphism:
\begin{equation}\label{frac:eq:inactive-product-splitting}
 X\simeq Z\times\prod_{j\in J_0}\PP^{m_j}.
\end{equation}
It remains only to record the properties of $Z$.  Since the second
factor in \eqref{frac:eq:inactive-product-splitting} is smooth, the
standard product formulas for canonical divisors and discrepancies show
that $Z$ is $\QQ$-Gorenstein and klt.  Restricting the anticanonical
decomposition to a slice $Z\times\{p\}$ gives
\[
 -K_Z\equiv
 \sum_{i\notin J_0}(m_i+1-\eps_i)G_i,
\]
where the $G_i$ are nef and their sum is ample.  Hence $Z$ is Fano.
The product decomposition also gives $\rho(Z)=\rho(X)-|J_0|$.
Numerical independence and optimality of the descended decomposition
follow immediately from those of the original decomposition: any
improvement on $Z$, together with the standard anticanonical
decompositions of the projective-space factors, would improve the
original decomposition on $X$.  Finally, Lemma~\ref{lem:product-index} gives $c_Z=c_X$.
This proves the lemma.
\end{proof}

\section{Complexity one}\label{sec:one}

Throughout this section $X$ is a klt Fano variety with $c_X=1$.  Put
$n:=\dim X$ and $\rho:=\rho(X)$, and retain the rounding notation
$m_i,\eps_i,M,e$ from Section~\ref{sec:optimal}.

\subsection{Integral coefficients and active directions}
By Corollary~\ref{frac:cor:reduction}, every optimal vertex decomposition
has $M=n-1$, $s=\rho(X)$, and $e=0$. Since the defects are nonnegative,
each $\eps_i=0$. Thus, fixing such a decomposition, we have
\begin{equation}\label{eq:integer}
 -K_X\equiv\sum_{i=1}^{\rho}(m_i+1)L_i,
 \qquad m_i\in\ZZ_{\geq0},\qquad\sum_i m_i=n-1.
\end{equation}
Set $m=(m_1,\ldots,m_\rho)$ and
$L^\alpha:=\prod_iL_i^{\alpha_i}$.  As in
Section~\ref{frac:sec:mixed}, define
\[
 d_i:=L^{m+\ee_i},
 \qquad
 I:=\{i\mid d_i>0\},
 \qquad
 J:=\{i\mid d_i=0\},
 \qquad
 J_0:=\{j\in J\mid\eps_j=0\}.
\]
The common identity \eqref{frac:eq:weighted} specializes to
\begin{equation}\label{eq:two-units}
 d_i\in\ZZ_{\geq0},\qquad\sum_i d_i=2.
\end{equation}
Hence the active set $I$ has one or two elements. Since $\eps_j=0$ for every
$j\in J$, we have $J_0=J$ in Lemma~\ref{frac:lem:inactive-split}. That lemma gives
\begin{equation}\label{eq:splitting}
 X\simeq Z\times\prod_{j\in J}\PP^{m_j},
 \qquad c_Z=1,\qquad\rho(Z)=|I|\leq2,
\end{equation}
where $Z$ is klt Fano and the active divisors descend to an optimal,
numerically independent decomposition on $Z$. Each inactive divisor is
the hyperplane class of its displayed factor. Since the product of the
top hyperplane powers on these factors has degree one, each active
$d_i$ is also the corresponding mixed intersection on $Z$.
It therefore remains to classify the factors of Picard number one and
two.

\subsection{Picard number one}
\begin{proposition}\label{prop:rho-one}
If $\rho(X)=1$ and $c_X=1$, then $X$ is a normal quadric of Picard
number one.
\end{proposition}
\begin{proof}
The curve case has $c_X=0$, so $n\geq2$. Equation~\eqref{eq:integer} gives
\[
 -K_X\equiv nH
\]
for an ample Cartier divisor $H$. By Lemma~\ref{lem:descent} this is a
$\QQ$-linear equivalence. Araujo--Druel
\cite[Theorem~1.1\ (2)]{AraujoDruel}, applied with zero boundary,
therefore identifies $(X,H)$ with a possibly singular quadric
$(Q^n,\OO_{Q^n}(1))$. Since $\rho(X)=1$, the smooth quadric surface is excluded.
\end{proof}

\subsection{Auxiliary geometric facts}
\begin{lemma}\label{lem:degree-one}
Let $X$ be a normal projective variety of dimension $n\geq2$ with
$K_X$ $\QQ$-Cartier, and let $H$ be a nef and big Cartier divisor. Suppose
\[
 H^n=1,
 \qquad
 g(X,H):=1+\frac12(K_X+(n-1)H)H^{n-1}=0.
\]
Then $|H|$ is base-point-free, $h^0(X,H)=n+1$, and its morphism
$X\to\PP^n$ is birational.
\end{lemma}
\begin{proof}
By \cite[Corollary~(2.11)]{Fujita89}, $\Delta(X,H)=0$.  Hence
\cite[Theorem~(1.1)]{Fujita89} shows that $H$ is globally generated and
defines a birational morphism.  Since $H^n=1$ and
$\Delta(X,H)=0$, we have $h^0(X,H)=n+1$, so its image is $\PP^n$.
\end{proof}

\begin{lemma}\label{lem:bilinear-singularities}
Let $a,b\geq2$, and let $Y\subset\PP^a\times\PP^b$ be an integral
hypersurface of bidegree $(1,1)$. Then $Y$ is normal, Gorenstein, and
terminal.
\end{lemma}
\begin{proof}
After independent coordinate changes the equation has the normal form
\begin{equation}\label{eq:bilinear}
 x_0y_0+\cdots+x_{r-1}y_{r-1}=0,
 \qquad 2\leq r\leq\min(a+1,b+1),
\end{equation}
and the Jacobian computation gives
$\Sing(Y)=\PP^{a-r}\times\PP^{b-r}$; see
\cite[Section~3.2 and Lemma~3.1]{SaitoSegre}.  Thus $Y$ is normal.  Being a
Cartier hypersurface in a smooth variety, it is Gorenstein by
\cite[Theorem~21.3]{EisenbudCA}.  Locally along $\Sing(Y)$ it is a smooth
product with $Q=\{\sum_{i=1}^ru_iv_i=0\}\subset\mathbb A^{2r}$; blowing up
the vertex gives
$K_{\widetilde Q}=\pi^*K_Q+(2r-3)E$, so $Q$, and hence $Y$, is terminal.
\end{proof}

\subsection{Picard number two}
\begin{proposition}\label{prop:rho-two}
Let $X$ be klt Fano of dimension $n\geq2$, and let $L,M$ be independent
nef Cartier divisors with
\begin{equation}\label{eq:rank-two}
 -K_X\equiv aL+bM,\quad a,b\in\ZZ_{>0},\quad a+b=n+1,
 \qquad L^{a-1}M^b=L^aM^{b-1}=1.
\end{equation}
Then, after interchanging $L,M$ if necessary, exactly one of the following
holds:
\begin{enumerate}[label=(\roman*)]
\item $a,b\geq2$ and $X$ is an integral $(1,1)$ hypersurface in
$\PP^a\times\PP^b$;
\item $X\simeq\Bl_{\PP^{a-2}}\PP^n$, with $L$ the hyperplane pullback
and $M=L-E$.
\end{enumerate}
\end{proposition}
\begin{proof}
Set \(A:=L+M\) and \(q_k:=L^kM^{n-k}\) for \(0\leq k\leq n\).
By Lemmas~\ref{lem:ample-sum} and \ref{lem:exchange}, $A$ is ample, the
$q_k$ are nonnegative integers, and their positive support is an interval
$[u,v]\cap\ZZ$.  Since $a+b=n+1$, \eqref{eq:rank-two} gives
$q_{a-1}=q_a=1$, so $u\leq a-1<a\leq v$.  Log-concavity and the two
consecutive values $q_{a-1}=q_a=1$ imply inductively that every positive
$q_k$ is at most $1$; hence
\begin{equation}\label{eq:ranktwo-ones}
 q_k=1\qquad(u\leq k\leq v).
\end{equation}

The divisors $L$ and $M$ cannot both be big.  Otherwise $u=0$, $v=n$, and
\eqref{eq:ranktwo-ones} gives $(tL+sM)^n=(t+s)^n$.  With $D=L-M$ this
implies $DA^{n-1}=D^2A^{n-2}=0$, so the numerical Hodge index theorem
\eqref{eq:hodge} gives $D\equiv0$, contrary to the independence of $L$ and
$M$.  After interchanging them if necessary, assume that $M$ is not big;
thus $u>0$.

Let $g:X\to T$ be the contraction defined by the semiample divisor $M$.
By Lemmas~\ref{lem:nef-monoid} and \ref{lem:descent}, there is an ample
Cartier divisor $H$ on $T$ such that
\begin{equation}\label{eq:M-pullback-ranktwo}
 M\sim g^*H.
\end{equation}
We determine $\dim T$ directly from the $q_k$.  By
\eqref{eq:numerical-dimension-definition},
\[
 \nu(M)=\max\{j\mid M^jA^{n-j}>0\},\qquad
 M^jA^{n-j}=\sum_{k=0}^{n-j}\binom{n-j}{k}q_k.
\]
Since all $q_k\geq0$ and their positive support is $[u,v]\cap\ZZ$, the last
sum is positive exactly when $u\leq n-j$.  Hence $\nu(M)=n-u$, and
Lemma~\ref{lem:numdim-basic} gives $\dim T=n-u$.  Thus a general fibre $F$
has
\begin{equation}\label{eq:ranktwo-fibre-dim}
 \dim F=u.
\end{equation}
For such $F$, adjunction and \eqref{eq:M-pullback-ranktwo} give
$-K_F\equiv aL|_F$, with $L|_F$ ample.  Lemma~\ref{lem:index} yields
$a\leq u+1$, whereas $q_{a-1}>0$ gives $u\leq a-1$.  Therefore equality
holds, and
\begin{equation}\label{eq:fibre-ranktwo}
 u=a-1\geq1,
 \qquad
 (F,L|_F)\simeq(\PP^{a-1},\OO_{\PP^{a-1}}(1)),
 \qquad
 \dim T=b.
\end{equation}
By the projection formula, $q_{a-1}=1$ and
$(L|_F)^{a-1}=1$ imply
\begin{equation}\label{eq:H-degree-ranktwo}
 H^b=1.
\end{equation}
Moreover, for $1\leq t\leq b-1$, the vector $(0,t)$ lies in the vanishing
box, so Lemma~\ref{lem:Bm-vanishing} gives
$\chi(X,\OO_X(-tM))=0$.  Lemma~\ref{lem:image-P}(ii) therefore
shows that
\begin{equation}\label{eq:base-ranktwo}
 (T,H)\simeq(\PP^b,\OO_{\PP^b}(1)).
\end{equation}
Thus $g:X\to\PP^b$ and $M=g^*\OO_{\PP^b}(1)$.

\smallskip\noindent
\emph{Case 1: $L$ is not big.}
Interchanging $L$ and $M$ in the preceding argument gives
$f:X\to\PP^a$ with $L=f^*\OO_{\PP^a}(1)$; in particular $a,b\geq2$.
The product morphism $\Phi=(f,g):X\to\PP^a\times\PP^b$ is finite onto its
image because $\Phi^*\OO(1,1)=L+M$ is ample.  Let $p_1$ and $p_2$ be the
projections of $\PP^a\times\PP^b$, and put
$H_1:=p_1^*\OO_{\PP^a}(1)$ and $H_2:=p_2^*\OO_{\PP^b}(1)$.  Its image
$W$ has dimension $n=a+b-1$, hence is an integral Cartier divisor, say
$[W]=rH_1+sH_2$.  If $d=\deg\Phi$, then the projection formula gives
\[
 1=L^{a-1}M^b=dr,
 \qquad
 1=L^aM^{b-1}=ds.
\]
Thus $d=r=s=1$.  Hence $W$ is an integral $(1,1)$ hypersurface and
$\Phi:X\to W$ is finite birational.  Lemma~\ref{lem:bilinear-singularities}
shows that $W$ is normal, so $\Phi$ is an isomorphism.  This is
(i).

\smallskip\noindent
\emph{Case 2: $L$ is big.}
Then $v=n$, so \eqref{eq:ranktwo-ones} gives
\begin{equation}\label{eq:L-degrees-ranktwo}
 L^n=1,
 \qquad
 L^{n-1}M=1.
\end{equation}
Using $a+b=n+1$, we obtain
\[
 g(X,L)=1+\frac12(K_X+(n-1)L)L^{n-1}
       =1+\frac12(n-1-a-b)=0.
\]
Lemma~\ref{lem:degree-one} therefore gives a birational morphism
\begin{equation}\label{eq:f-ranktwo}
 f:X\longrightarrow\PP^n,
 \qquad
 L=f^*\OO_{\PP^n}(1).
\end{equation}

Set $\mathcal F:=f_*\OO_X(M)$.  Since $f$ is birational and $\PP^n$ is
smooth, $\mathcal F$ is rank-one torsion-free and
$\mathcal F^{**}\simeq\OO_{\PP^n}(k)$ for some $k\in\ZZ$.  On a general
line contained in the locus where $f$ is an isomorphism and
$\mathcal F=\mathcal F^{**}$, \eqref{eq:L-degrees-ranktwo} gives
$k=ML^{n-1}=1$.  Hence $\mathcal F^{**}\simeq\OO_{\PP^n}(1)$, and
\begin{equation}\label{eq:linear-system-injection}
 H^0(X,M)=H^0(\PP^n,\mathcal F)
 \hookrightarrow H^0(\PP^n,\OO_{\PP^n}(1)).
\end{equation}
By \eqref{eq:base-ranktwo}, the left-hand side has dimension $b+1$.
Thus it defines a linear projection
$\pi_V:\PP^n\dashrightarrow\PP^b$ whose centre is a linear subspace
$\Lambda\simeq\PP^{n-b-1}=\PP^{a-2}$, and on the common domain of
definition $g\circ f^{-1}=\pi_V$.

Again $\Phi=(f,g):X\to\PP^n\times\PP^b$ is finite onto its image because
$L+M$ is ample.  Its image is the closure of the graph of $\pi_V$, namely
$\Bl_\Lambda\PP^n$.  Hence $\Phi$ is a finite birational morphism onto the
smooth variety $\Bl_\Lambda\PP^n$, and therefore an isomorphism.  If $E$
is the exceptional divisor, the resolved projection is given by the strict
transforms of the hyperplanes containing $\Lambda$, so $M\sim L-E$.
Thus
\[
 X\simeq\Bl_{\PP^{a-2}}\PP^n,
\]
with $L$ the hyperplane pullback and $M=L-E$.  This is (ii), and
the two cases exhaust the possibilities.
\end{proof}

\section{A birational degree lemma}\label{frac:sec:birational}

This section is independent of the standing optimal-decomposition
notation. The letters $A,B,\eps,b,\delta$, and $E$ below are local to the
birational morphism in the lemma and do not redefine the symbols used in
the surrounding classification sections.

\begin{lemma}\label{frac:lem:birational}
Let $X$ be a complex klt Fano $n$-fold, $n\geq2$, and let
$f:X\to\PP^n$ be a birational morphism. Let $H$ be a hyperplane
divisor on $\PP^n$, and put $A=f^*H$.
Suppose a nonzero nef Cartier divisor $B$ on $X$ and real numbers
$b,\eps$ satisfy
\[
 -K_X\equiv(n-\eps)A+bB,\qquad
 \delta=BA^{n-1}\in\ZZ_{\geq2},\qquad
 \delta b=1+\eps,\qquad0\leq\eps<b.
\]
Then either $f$ is an isomorphism and $B\sim\delta A$, or $f$ is the
ordinary blow-up of a codimension-two linear subspace and
\[
 \eps=0,\qquad b=\frac1\delta,\qquad B\sim\delta(A-E),
\]
where $E$ is its exceptional Cartier divisor. In the latter case, replacing $bB$ by $\delta b(A-E)$ strictly increases the coefficient sum
in a nef Cartier decomposition.
\end{lemma}
\begin{proof}
The hypotheses give $b>0$ and $(\delta-1)b<1$. Since $\delta\geq2$,
we have $0\leq\eps<b<1$ and $n-\eps>0$. Thus, with $\lambda=1/b$,
\begin{equation}\label{frac:eq:lambda-lower}
 \lambda>\delta-1,\qquad 2\lambda>\delta,\qquad 0\leq\eps<1.
\end{equation}
Moreover $A+B$ is ample by Lemma~\ref{lem:ample-sum}, and $B$ is
$f$-ample because $-K_X\equiv_f bB$. Purity of the exceptional locus
\cite[Chapter~2, Corollary~2.63]{KM} applies because $X$ is projective
and $\PP^n$ is smooth, hence $\QQ$-factorial. In particular, if $f$
has no exceptional prime, it is an isomorphism. Then
$\Pic(X)=\ZZ[A]$ and $BA^{n-1}=\delta$ give $B\sim\delta A$.

Assume that $f$ has an exceptional prime. Choose canonical divisors
on $X$ and $\PP^n$ using the same nonzero rational top form, and let
$\{E_\alpha\}$ be all the $f$-exceptional prime divisors on $X$.
Write $R:=K_X-f^*K_{\PP^n}=\sum_\alpha a_\alpha E_\alpha$.
Since $\PP^n$ is smooth, every $a_\alpha=a(E_\alpha,\PP^n)$ is a
positive integer: at the generic point of $E_\alpha$, where $X$ is
smooth, it is the positive vanishing order of the Jacobian determinant.
Thus $\Supp R$ contains every exceptional prime, and purity gives
$\operatorname{Exc}(f)=\Supp R$.
Choose a Cartier representative of $B$, still denoted by $B$, and put
$G=f_*B$. The divisor $G$ is Cartier because $\PP^n$ is smooth.
By the projection formula, its degree is $BA^{n-1}=\delta$, so
$G\sim\delta H$. The Cartier divisor $B-f^*G$ is exceptional and is
unchanged if the representative of $B$ is changed by a principal divisor.
The canonical-class identity gives $R\equiv b(\delta A-B)$.
Consequently $B-f^*G+\lambda R$ is exceptional and numerically trivial
over $\PP^n$. Applying the negativity lemma
\cite[Lemma~4.16]{FujinoFundamental} to it and its negative gives zero.
It follows that
\begin{equation}\label{frac:eq:birational-divisor}
 B\sim\delta A-\lambda R,
 \qquad \lambda R\text{ is an effective integral Cartier divisor},
\end{equation}
where integrality and the Cartier property follow from
$\lambda R=f^*G-B$.
Since $f_*\OO_X=\OO_{\PP^n}$ and $\lambda R$ is effective,
$\mathcal I:=f_*\OO_X(-\lambda R)$ is an ideal subsheaf of
$\OO_{\PP^n}$. The projection formula gives
$\mathcal F:=f_*\OO_X(B)\simeq\mathcal I(\delta)$.
For an ideal sheaf, $\ord_{E_\alpha}$ means the minimum vanishing
order of the pullbacks of its local generators at the generic point
of the centre of $E_\alpha$. The definition of $\mathcal I$ gives
\begin{equation}\label{frac:eq:ideal-valuations}
 \ord_{E_\alpha}(\mathcal I)\geq\lambda a_\alpha.
\end{equation}
Set $Z=\bigcup_\alpha f(E_\alpha)=f(\operatorname{Exc}(f))$.
The positive orders in \eqref{frac:eq:ideal-valuations} force every
local section of $\mathcal I$ to vanish along every centre. Over
$\PP^n\setminus Z$, the morphism $f$ is an isomorphism and
$\mathcal I=\OO_{\PP^n}$. Thus
$\Supp(\OO_{\PP^n}/\mathcal I)=Z$, a nonempty closed set of
codimension at least two.

\smallskip\noindent
\emph{Global generation.}
Since $B-K_X\equiv_f(1+b)B$ is $f$-ample, relative
Kawamata--Viehweg vanishing gives $R^qf_*\OO_X(B)=0$ for $q>0$.
The projection formula and Leray therefore give
\[
 H^i(\PP^n,\mathcal F(-i))
 \simeq H^i(X,\OO_X(B-iA)).
\]
For $1\leq i<n$, the divisor
$(B-iA)-K_X\equiv(n-i-\eps)A+(1+b)B$ is ample: both coefficients
are positive, $A,B$ are nef, and $A+B$ is ample. Kawamata--Viehweg
vanishing therefore makes these groups zero.
Since klt singularities are Cohen--Macaulay, Serre duality for $i=n$
identifies the dual group with
$H^0(X,\omega_X\otimes\OO_X(nA-B))$.
A nonzero section would give an effective Weil divisor linearly
equivalent to $K_X+nA-B$. However, $A^n=1$ by birationality, and
\[
 (K_X+nA-B)A^{n-1}
 =\eps-(1+b)\delta=-\delta-1<0,
\]
contradicting the nonnegative intersection of an effective divisor
with the nef class $A^{n-1}$. Hence
$\mathcal F=\mathcal I(\delta)$ is $0$-regular and globally generated
\cite[Proposition~3.8~(3)]{BallicoMalaspina}.
Its sections define a degree-$\delta$ system $\mathcal M$ and generate
the ideal $\mathcal I$ after local trivialization of $\OO(\delta)$.
Thus the base ideal of $\mathcal M$ is exactly $\mathcal I$, its
set-theoretic base locus is $Z$, and it is movable. Moreover,
$\ord_{E_\alpha}(\mathcal M)=\ord_{E_\alpha}(\mathcal I)$,
where the left side is the minimum coefficient of $E_\alpha$ in
the pullbacks of members of $\mathcal M$.

\smallskip\noindent
\emph{The linear centre.}
By \eqref{frac:eq:ideal-valuations} and $\lambda>\delta-1$, the
discrepancy of $E_\alpha$ for $(\PP^n,(\delta-1)^{-1}\mathcal M)$ is
\[
 a_\alpha-\frac{\ord_{E_\alpha}(\mathcal M)}{\delta-1}
 \leq a_\alpha\left(1-\frac{\lambda}{\delta-1}\right)<0.
\]
Here $\operatorname{mult}_z\mathcal M$ denotes the minimum
multiplicity at $z$ of its members. If this were at most $\delta-1$
at a point $z\in Z$, upper semicontinuity would give the same bound
at every point of an open neighbourhood of $z$. The multiplicity
estimate \cite[Lemma~28]{KollarRigidity}, applied to that neighbourhood
with coefficient $(\delta-1)^{-1}$, would make the restricted pair
canonical, contradicting the negative discrepancy of an $E_\alpha$
whose centre contains $z$. Hence
$\operatorname{mult}_z\mathcal M>\delta-1$ for every $z\in Z$.
Multiplicities are integers, and a nonzero degree-$\delta$ form has
multiplicity at most $\delta$ at any point. Therefore every member
has multiplicity exactly $\delta$ at every point of $Z$.

The multiplicity-$\delta$ locus of a degree-$\delta$ form is its
linear vertex: it is cut out by all partial derivatives of order
$\delta-1$, which are linear forms. The lower derivatives vanish
there as well by Euler's identity. The intersection of the vertices
of all members contains $Z$ and is contained in their common zero
locus, which is $Z$. Thus $Z$ is a linear subspace. Give it its
reduced scheme structure, write $C=Z\simeq\PP^{n-k}$ with $k\geq2$,
and let $p:S=\Bl_C\PP^n\to\PP^n$ have exceptional divisor $E_0$.
In homogeneous coordinates with $C=(x_1=\cdots=x_k=0)$, every
defining form depends only on the normal variables $x_1,\ldots,x_k$:
multiplicity $\delta$ along $C$ excludes all monomials of smaller
total normal degree. These degree-$\delta$ forms have no common zero
in $\PP^{k-1}$, since such a zero would give a base point outside $C$.
After pulling back to $S$ and removing the fixed divisor $\delta E_0$,
the residual system is therefore the pullback of a base-point-free
system on $\PP^{k-1}$ under the resolved linear projection from $C$.
Since $\mathcal M$ generates $\mathcal I$, this proves the equality
of ideal sheaves
\begin{equation}\label{frac:eq:ideal-on-blowup}
 \mathcal I\OO_S=\OO_S(-\delta E_0).
\end{equation}

\smallskip\noindent
\emph{The exceptional valuation and the blow-up.}
For an exceptional prime $E$ of $f$, put $v=\ord_E$ and
$a=a(E,\PP^n)$. Evaluate $v$ on divisors on $S$ using a common
resolution of $X$ and $S$ over $\PP^n$. Since the centre of $v$ on
$\PP^n$ lies in $C$ and $\mathcal I_C\OO_S=\OO_S(-E_0)$, we have
$w:=v(\mathcal I_C)=v(E_0)>0$.
Equations~\eqref{frac:eq:ideal-valuations} and
\eqref{frac:eq:ideal-on-blowup} give
\begin{equation}\label{frac:eq:valuation-bound}
 \lambda a\leq v(\mathcal I)=\delta w.
\end{equation}
If $v\ne\ord_{E_0}$, its centre on $S$ lies in $E_0$ and has
codimension at least two: a normalized divisorial valuation with a
codimension-one centre on the normal variety $S$ is the order along
that prime divisor. Thus $E$ is exceptional over $S$. Since $S$ and
$E_0$ are smooth, \cite[Chapter~2, Corollary~2.31~(3)]{KM} gives
$a(E,S)-v(E_0)=a(E,S,E_0)\geq0$.
Using $K_S=p^*K_{\PP^n}+(k-1)E_0$ on the common resolution yields
\[
 a=a(E,S)+(k-1)w\geq kw\geq2w,
\]
contradicting \eqref{frac:eq:valuation-bound} and $2\lambda>\delta$.
Thus every exceptional prime defines the valuation $\ord_{E_0}$.
Distinct primes on the normal variety $X$ define distinct normalized
valuations, so there is exactly one, denoted $E$. For this valuation
$w=1$ and $a=k-1$; the inequality
$\lambda(k-1)\leq\delta<2\lambda$ forces $k=2$.
Consequently $R=E$. Since $\lambda R$ is integral, $\lambda$ is an
integer, and $\delta-1<\lambda\leq\delta$ gives $\lambda=\delta$.
Therefore
\[
 b=\frac1\delta,\qquad \eps=0,\qquad B\sim\delta(A-E).
\]

The divisor $-\delta E$ is Cartier by
\eqref{frac:eq:birational-divisor}, and it is $f$-ample because it
is linearly equivalent to $B-\delta A$.
Normality and the uniqueness of the exceptional valuation give, for
every integer $m\geq0$,
\[
 f_*\OO_X(-m\delta E)=\mathcal I_C^{m\delta}.
\]
Indeed, on an open set $U\subset\PP^n$ meeting $C$, sections of
the left-hand side are precisely the regular functions
$h\in\OO_{\PP^n}(U)$ satisfying $\ord_{E_0}(h)\geq m\delta$.
Regular functions downstairs pull back regularly, and the only
additional divisorial vanishing condition is the one along $E$;
conversely, $f_*\OO_X=\OO_{\PP^n}$ ensures regularity downstairs.
Since $C$ is a smooth linear subspace, $\ord_{E_0}$ measures total
degree in local coordinates normal to $C$. The condition excludes
all terms of degree less than $m\delta$ and hence gives the ordinary
ideal power. Away from $C$, both sheaves are the structure sheaf.
The $f$-ample line bundle $\OO_X(-\delta E)$ reconstructs $X$ as
the relative Proj of its section algebra. Hence
\[
 X\simeq\operatorname{Proj}_{\PP^n}
 \bigoplus_{m\geq0}\mathcal I_C^{m\delta}
 \simeq\Bl_C\PP^n,
\]
since this algebra is the $\delta$th Veronese of the Rees algebra
$\bigoplus_{j\geq0}\mathcal I_C^j$, and passage to a Veronese
subalgebra does not change Proj. In particular $E$ is Cartier and
$A-E$ is the nonzero nef Cartier divisor defining the resolved
linear projection from $C$.
Replacing $bB$ by $\delta b(A-E)$ increases the coefficient sum by
$(\delta-1)b>0$, as claimed.
\end{proof}

\section{The active factors and the classification}\label{frac:sec:classification}

We retain from Section~\ref{frac:sec:mixed} the optimal decomposition
\eqref{frac:eq:mixed-data} and the notation
\[
 m=(m_1,\ldots,m_s),
 \qquad
 e=\sum_i\eps_i,
 \qquad
 L^\alpha:=\prod_iL_i^{\alpha_i},
 \qquad
 d_i:=L^{m+\ee_i},
\]
\[
 I:=\{i\mid d_i>0\},
 \qquad
 J:=\{i\mid d_i=0\}.
\]
Assume throughout this section that $1<c_X<2$. Thus
$0<e=c_X-1<1$.

\subsection{The possible active degrees}

\begin{lemma}\label{frac:lem:degrees}
The active degrees are either a single degree $d\geq2$, or two degrees
$(d,1)$ with $d\geq1$. The possibility $(1,1)$ does not occur when
$e>0$.
\end{lemma}
\begin{proof}
Put
\[
 D:=\sum_i d_i,
 \qquad
 d:=\max_i d_i.
\]
Only active indices contribute to these sums.  Expanding
\eqref{frac:eq:weighted} gives
\(D-\sum_i\eps_i d_i=2\), equivalently
\(D-2=\sum_i\eps_i d_i\).  Since $d_i\leq d$ for every $i$ and
$e=\sum_i\eps_i<1$,
\[
 D-2\leq d\sum_i\eps_i=de<d.
\]
Hence $D-d<2$.  The number $D-d$ is a nonnegative integer and is the sum of all active
degrees after removing one largest degree.  Therefore it is either zero
or one.  If it is zero, there is exactly one active index.  If it is one,
there are exactly two active indices, and the smaller active degree is
one.  Thus the only numerical possibilities are a single degree $d$ or
two degrees $(d,1)$.

A single active degree cannot equal one.  Indeed, in that case
\eqref{frac:eq:weighted} would read
\(1-\eps_p=2,\)
which is impossible because $0\leq\eps_p<1$.  Hence in the one-active
case $d\geq2$.

It remains to exclude the two-active possibility $(1,1)$.  Suppose
$d_p=d_q=1$.  Equation~\eqref{frac:eq:weighted} becomes
\((1-\eps_p)+(1-\eps_q)=2.\)
Since the defects are nonnegative,
\begin{equation}\label{frac:eq:two-unit-zero-defect}
 \eps_p=\eps_q=0.
\end{equation}
On the other hand $e>0$ because $1<c_X<2$.  Hence there exists an
inactive index $j\in J$ with
\(\eps_j>0.\)
By Lemma~\ref{frac:lem:inactive-base},
\(m_j>0\) and \(\nu(L_j)=m_j\).  Consider the $(n-1)$-fold intersection class
\[
 C:=L^{m-\ee_j+\ee_p}.
\]
It is well defined because $m_j>0$.  We have
\(L_jC=L^{m+\ee_p}=d_p=1.\)
If $k\in J$ and $k\ne j$, then
\(L_kC=L^{m-\ee_j+\ee_p+\ee_k}=0.\)
To see this, suppose the intersection were positive.  In the expansion
of
\[
 L_k^{m_k+1}\left(\sum_iL_i\right)^{n-m_k-1}
\]
all terms are nonnegative and one of them is a positive multiple of
$L_kC$.  This would give
$\nu(L_k)\geq m_k+1$, contradicting
Lemma~\ref{frac:lem:inactive-base}.

Intersect the anticanonical decomposition with $C$.  By
\eqref{frac:eq:two-unit-zero-defect}, the coefficients of the two active
classes are integers.  Every inactive class other than $L_j$ has zero
intersection with $C$, while $L_jC=1$.  Thus
\[
 -K_XC=N-\eps_j
\]
for some integer $N$.  However \eqref{eq:K-integral}, applied to
$u=m-\ee_j+\ee_p$, gives
\[
 K_XC\in\ZZ.
\]
It follows that $\eps_j\in\ZZ$, contradicting
$0<\eps_j<1$.  Hence $(1,1)$ cannot occur.

Consequently the active degrees are either a single degree $d\geq2$, or
two degrees $(d,1)$ with the case $d=1$ excluded by the preceding
argument.
\end{proof}

\subsection{A surface obstruction}

\begin{lemma}\label{frac:lem:surface}
Let $S$ be a klt Fano surface.
\begin{enumerate}[label=(\roman*)]
\item If an extremal ray $R$ defines a birational contraction and
$R=\RR_{\geq0}[C]$ for an integral curve $C$, then
$0<-K_S\cdot C\leq1$.
\item There do not exist nef Cartier divisors $A_1,\ldots,A_q,C$,
where $q\geq1$, and rational numbers $a_1,\ldots,a_q,\eta$ satisfying
\[
 -K_S\equiv\sum_{\ell=1}^q a_\ell A_\ell+(2-\eta)C,
 \qquad C\not\equiv0,\quad C^2=0,
 \qquad 0<\eta<\min\{1,a_1,\ldots,a_q\}.
\]
\end{enumerate}
\end{lemma}
\begin{proof}
Part (i) follows from \cite[Lemma~4.8]{DeVleming}; for a surface, the
fibres of a birational contraction have dimension at most one.

For (ii), suppose such a decomposition exists. Since $C\not\equiv0$ is
nef and $C^2=0$, we have $\rho(S)\geq2$. Choose an extremal ray $R$
with $C\cdot R>0$, and let $\Gamma$ be an integral curve generating
$R$. Since $C$ is Cartier, $C\Gamma\geq1$, and hence
\[
 -K_S\Gamma
 =\sum_{\ell=1}^q a_\ell A_\ell\Gamma+(2-\eta)C\Gamma
 \geq2-\eta>1.
\]
By (i), the contraction of $R$ is not birational. It cannot contract
$S$ to a point since $\rho(S)\geq2$, so it is a contraction onto a
curve. Let $F$ be a general fibre. Then $F\simeq\PP^1$ and
$-K_SF=2$. Intersecting the assumed decomposition with $F$ gives
\[
 2=\sum_{\ell=1}^q a_\ell A_\ell F+(2-\eta)CF.
\]
Since $CF$ is a positive integer and $2-\eta>1$, we must have $CF=1$,
and hence
\[
 \sum_{\ell=1}^q a_\ell A_\ell F=\eta.
\]
But every $A_\ell F$ is a nonnegative integer. Thus the left-hand side
is either zero or at least $\min_\ell a_\ell>\eta$, a contradiction.
\end{proof}

We use general product lines in the inactive base and cut their inverse
images by general members of base-point-free systems.  We also use
relative spannedness in the form of \cite[Theorem~1.1]{FujinoSpanned}.

\subsection{One active degree}

\begin{proposition}\label{frac:prop:one-active}
Suppose $I=\{p\}$ and $d_p=d$. Then $\eps_j=0$ for every $j\in J$,
$d\geq3$, and
\[
 X\simeq Z\times\prod_{j\in J}\PP^{m_j},\qquad\rho(Z)=1.
\]
Writing $r=m_p+1$, the primitive ample Cartier class $H_Z$ on $Z$ satisfies
\begin{equation}\label{frac:eq:rank-one}
 \dim Z=r\geq2,
 \qquad
 -K_Z\equiv\left(r-1+\frac2d\right)H_Z,
 \qquad H_Z^r=d.
\end{equation}
\end{proposition}
\begin{proof}
Let $I=\{p\}$ and write $d_p=d$.  Since $d_j=0$ for every $j\in J$,
\eqref{frac:eq:weighted} becomes
\((1-\eps_p)d=2.\)
Hence
\begin{equation}\label{frac:eq:one-active-eps}
 \eps_p=1-\frac2d.
\end{equation}
Because
$e=\eps_p+\sum_{j\in J}\eps_j<1$, we also have
\begin{equation}\label{frac:eq:one-active-inactive-defect-bound}
 \sum_{j\in J}\eps_j<\frac2d.
\end{equation}

We first prove that $\eps_j=0$ for every $j\in J$.  Suppose to the
contrary that some $j\in J$ has
\(\eta:=\eps_j>0.\)
Then \eqref{frac:eq:one-active-inactive-defect-bound} gives
\begin{equation}\label{frac:eq:eta-one-active}
 0<\eta<\min\left\{1,\frac2d\right\}.
\end{equation}
By Lemma~\ref{frac:lem:inactive-base}, the inactive contraction is
\[
 f:X\longrightarrow
 Y=\prod_{k\in J}\PP^{m_k}.
\]
Choose a general line $\ell\subset\PP^{m_j}$ in the $j$th factor and
general points in all other factors.  Let $X'$ be the inverse image of
the resulting copy of $\PP^1$ in $Y$.  Equivalently, $X'$ is obtained
by intersecting $X$ with $m_k$ general pullbacks of hyperplanes from
each factor $k\ne j$ and with $m_j-1$ general pullbacks from the $j$th
factor.  Since $f$ has connected fibres and the chosen base curve is
connected, $X'$ is connected.  Successive general members of these
base-point-free systems are normal and klt by
\cite[Lemma~5.17(1)]{KM}; since they are Cartier, their differents are
zero, and inversion of adjunction
\cite[Theorem~4.9(1)]{KollarSingMMP} shows that the corresponding pairs
are plt.  Thus repeated Bertini and adjunction show that $X'$ is normal
and klt.  Since a connected normal Noetherian scheme has disjoint irreducible
components, $X'$ is irreducible.

Put
\[
 r:=m_p+1,
 \qquad
 H_{X'}:=L_p|_{X'},
 \qquad
 C_{X'}:=L_j|_{X'}.
\]
The relative dimension of $f$ is
\(n-\sum_{k\in J}m_k=m_p+1=r,\)
so $\dim X'=r+1$.  All inactive divisor classes other than $L_j$
become trivial on the chosen base curve.  Adjunction subtracts one copy
of each hyperplane used in the cut.  Hence
\begin{equation}\label{frac:eq:one-active-line}
 -K_{X'}\equiv
 \left(r-1+\frac2d\right)H_{X'}+(2-\eta)C_{X'},
 \qquad
 C_{X'}=f'^*\OO_{\PP^1}(1),
\end{equation}
where $f':X'\to\PP^1$ is the induced morphism.  The intersection number
$d_p=L^{m+\ee_p}$ becomes
\begin{equation}\label{frac:eq:one-active-line-degree}
 H_{X'}^rC_{X'}=d.
\end{equation}
The divisor $H_{X'}+C_{X'}$ is ample, being the restriction of a positive
combination of the original supporting divisors.  Both coefficients in
\eqref{frac:eq:one-active-line} are positive by
\eqref{frac:eq:eta-one-active}; since $H_{X'}$ and $C_{X'}$ are nef, it follows
that $-K_{X'}$ is ample.  Thus $X'$ is a klt Fano variety.  Since $C_{X'}$ is
pulled back from the base, $H_{X'}$ is $f'$-ample.

We next show that $H_{X'}$ is globally generated.  Since $X'$ is integral and
$\PP^1$ is smooth, the dominant morphism $f'$ is flat; hence every fibre
has dimension $r$.  Put
\[
 \lambda:=r-1+\frac2d.
\]
Then $K_{X'}+\lambda H_{X'}\equiv_{f'}0$ and
\[
 r<\lambda+1.
\]
Fujino's relative spannedness theorem
\cite[Theorem~1.1]{FujinoSpanned} therefore applies and shows that $H_{X'}$
is generated by relative global sections over $\PP^1$.  Relative
Kawamata--Viehweg vanishing gives
\(R^if'_*\OO_{X'}(H_{X'})=0\) for \(i>0\),
because $H_{X'}-K_{X'}$ is $f'$-ample.  Moreover
\[
 H_{X'}-C_{X'}-K_{X'}
 \equiv
 \left(r+\frac2d\right)H_{X'}+(1-\eta)C_{X'}
\]
is ample by \eqref{frac:eq:eta-one-active}.  Hence
\(H^1(X',\OO_{X'}(H_{X'}-C_{X'}))=0.\)
By projection formula and Leray,
\[
 H^1\left(\PP^1,
 f'_*\OO_{X'}(H_{X'})\otimes\OO_{\PP^1}(-1)\right)=0.
\]
Thus the sheaf $f'_*\OO_{X'}(H_{X'})$ is $0$-regular on $\PP^1$ and is
globally generated.  Pulling back its evaluation map gives a surjection
\[
 H^0\!\left(\PP^1,f'_*\OO_{X'}(H_{X'})\right)\otimes\OO_{X'}
 \longrightarrow f'^*f'_*\OO_{X'}(H_{X'}).
\]
By relative spannedness, the relative evaluation map
\[
 f'^*f'_*\OO_{X'}(H_{X'})\longrightarrow\OO_{X'}(H_{X'})
\]
is also surjective.  Composing these two maps and using
\[
 H^0\!\left(\PP^1,f'_*\OO_{X'}(H_{X'})\right)
 =
 H^0\!\left(X',\OO_{X'}(H_{X'})\right),
\]
we conclude that $H_{X'}$ is globally generated.

Choose $r-1$ general members of $|H_{X'}|$ and let $S$ be their complete
intersection; when $r=1$ there is no cut.  As above, repeated Bertini
and adjunction show that $S$ is a klt Fano surface satisfying
\begin{equation}\label{frac:eq:one-active-surface}
 -K_S\equiv\frac2dH_S+(2-\eta)C_S,
 \qquad
 H_S:=H_{X'}|_S,\quad C_S:=C_{X'}|_S.
\end{equation}
Since $C_{X'}$ comes from the base curve, $C_S^2=0$, while
\eqref{frac:eq:one-active-line-degree} gives $H_SC_S=d$.  In particular
$C_S\not\equiv0$, and $H_S,C_S$ are nef Cartier divisors.
By \eqref{frac:eq:eta-one-active}, Lemma~\ref{frac:lem:surface}
(ii) applies with $q=1$, $A_1=H_S$, and $a_1=2/d$, contradicting
\eqref{frac:eq:one-active-surface}. Therefore
\begin{equation}\label{frac:eq:all-inactive-integral-one-active}
 \eps_j=0
 \qquad(j\in J).
\end{equation}

Lemma~\ref{frac:lem:inactive-split} now splits off all inactive factors:
\[
 X\simeq Z\times\prod_{j\in J}\PP^{m_j},
 \qquad
 \rho(Z)=1,
 \qquad
 c_Z=c_X.
\]
Let $H_Z$ be the ample Cartier divisor on $Z$ to which $L_p$ descends;
equivalently, under the above product decomposition, $L_p$ is the pullback
of $H_Z$ from $Z$.  Since
\(\dim Z=n-\sum_{j\in J}m_j=m_p+1=r,\)
\eqref{frac:eq:one-active-eps} gives
\[
 -K_Z\equiv
 \left(m_p+1-\eps_p\right)H_Z
 =\left(r-1+\frac2d\right)H_Z.
\]
Since
\[
 d=d_p=L^{m+\ee_p}
   =L_p^r\prod_{j\in J}L_j^{m_j},
\]
and the top hyperplane degree on each factor
$\PP^{m_j}$ is one, the product decomposition gives
\[
 d=H_Z^r.
\]
If $H_Z$ were numerically divisible by an integer greater than one in the
Cartier lattice, replacing $H_Z$ by its primitive generator would enlarge
the coefficient in the anticanonical decomposition and contradict
optimality.  Hence $H_Z$ is primitive.

Finally $e>0$ and \eqref{frac:eq:all-inactive-integral-one-active} give
\(e=\eps_p=1-\frac2d>0,\)
so $d\geq3$.  The case $r=1$ is impossible, because a normal Fano curve
is $\PP^1$ and has nef complexity zero, whereas $c_Z=c_X>1$.  Thus
$r\geq2$, and all assertions follow.
\end{proof}

\subsection{Two active degrees}

\begin{proposition}\label{frac:prop:two-active}
There cannot be two active indices.
\end{proposition}
\begin{proof}
Assume that there are two active indices.  By
Lemma~\ref{frac:lem:degrees}, after relabelling we may write
\(d_p=d\geq2\) and \(d_q=1\).
Set
\(t:=\eps_q\) and \(b:=1-\eps_p\).
Equation~\eqref{frac:eq:weighted} becomes
\(db+(1-t)=2,\)
so
\begin{equation}\label{frac:eq:two-numbers}
 db=1+t.
\end{equation}
Furthermore
\[
 e=(1-b)+t+\sum_{j\in J}\eps_j<1,
\]
which is equivalent to
\begin{equation}\label{frac:eq:two-inactive-defect-bound}
 \sum_{j\in J}\eps_j<b-t.
\end{equation}
In particular $t<b$.  Since $d\geq2$ and $t<1$,
\[
 b=\frac{1+t}{d}<1.
\]
Thus
\begin{equation}\label{frac:eq:two-number-range}
 0\leq t<b<1.
\end{equation}

\smallskip\noindent
\emph{Step 1: the larger active degree has $m_p=0$.}
Suppose $m_p>0$.  Apply the mixed Khovanskii--Teissier inequality
\eqref{eq:mixed-KT} with
\[
 B=L_p,
 \qquad
 C=L_q,
 \qquad
 N_1\cdots N_{n-2}=L^{m-\ee_p}.
\]
This gives
\[
 \left(L^{m+\ee_q}\right)^2
 \geq
 L^{m+\ee_p}\,L^{m-\ee_p+2\ee_q}.
\]
Since the first two displayed active degrees are $1$ and $d$, we obtain
\(1\geq d\,L^{m-\ee_p+2\ee_q}.\)
The last intersection number is a nonnegative integer and $d\geq2$;
hence
\begin{equation}\label{frac:eq:two-KT-zero}
 L^{m-\ee_p+2\ee_q}=0.
\end{equation}

Put
\(u:=m-\ee_p+\ee_q.\)
Then $|u|=n-1$, so the parity relation \eqref{eq:parity} applies.
Using
\[
 K_X\equiv-\sum_i(m_i+1-\eps_i)L_i,
\]
we can write
\[
 K_X+\sum_i u_iL_i
 \equiv
 \sum_i\bigl(u_i-(m_i+1-\eps_i)\bigr)L_i.
\]
For the index $p$, since $u_p=m_p-1$, the coefficient of $L_p$ in
this expression is
\[
 (m_p-1)-(m_p+1-\eps_p)=\eps_p-2,
\]
and
\(L_pL^u=L^{m+\ee_q}=1.\)
For the index $q$, since $u_q=m_q+1$, the coefficient of $L_q$ is
\[
 (m_q+1)-(m_q+1-\eps_q)=\eps_q=t.
\]
However,
\(L_qL^u=L^{m-\ee_p+2\ee_q}=0\)
by \eqref{frac:eq:two-KT-zero}.  If $j\in J$, then $u_j=m_j$, so
\[
 L_jL^u
 =
 L_j^{m_j+1}\prod_{i\ne j}L_i^{u_i}.
\]
Hence $L_jL^u>0$ would imply
\[
 L_j^{m_j+1}\left(\sum_iL_i\right)^{n-m_j-1}>0,
\]
and therefore $\nu(L_j)\geq m_j+1$, contradicting
$\nu(L_j)=m_j$ from Lemma~\ref{frac:lem:inactive-base}.  Thus
$L_jL^u=0$.  Therefore
\[
 \left(K_X+\sum_i u_iL_i\right)L^u=\eps_p-2.
\]
The left-hand side belongs to $2\ZZ$ by \eqref{eq:parity}.  The right
hand side lies strictly between $-2$ and $-1$, because
$0<\eps_p=1-b<1$.  This contradiction proves
\begin{equation}\label{frac:eq:mp-zero}
 m_p=0.
\end{equation}

\smallskip\noindent
\emph{Step 2: $\eps_j=0$ for every $j\in J$.}
Suppose that an inactive index $j$ has
\(\eta:=\eps_j>0\).
By \eqref{frac:eq:two-inactive-defect-bound},
\begin{equation}\label{frac:eq:eta-two-active}
 0<\eta<b-t\leq\min\{b,1-t\}.
\end{equation}
Put \(r:=m_q+1\).  As in the proof of
Proposition~\ref{frac:prop:one-active}, use
Lemma~\ref{frac:lem:inactive-base} to restrict the inactive base to a
general line in the \(j\)-th factor and to general points in all other
factors.  Let
\[
 f':X'\longrightarrow\PP^1
\]
be the resulting morphism.  Since \(m_p=0\) by
\eqref{frac:eq:mp-zero} and \(|m|=n-1\), the relative dimension of the
inactive contraction is \(r\), and hence \(\dim X'=r+1\).  As in the
one-active case, repeated Bertini and adjunction show that \(X'\) is
normal and klt.

Put
\[
 A_{X'}:=L_q|_{X'},
 \qquad
 B_{X'}:=L_p|_{X'},
 \qquad
 C_{X'}:=L_j|_{X'}.
\]
The same adjunction calculation as in
Proposition~\ref{frac:prop:one-active} gives
\begin{equation}\label{frac:eq:two-line}
 -K_{X'}
 \equiv
 (r-t)A_{X'}+bB_{X'}+(2-\eta)C_{X'},
 \qquad
 C_{X'}=f'^*\OO_{\PP^1}(1).
\end{equation}
The divisors \(A_{X'},B_{X'},C_{X'}\) are nef Cartier, and their sum is
ample.  Since all coefficients in \eqref{frac:eq:two-line} are positive,
\(-K_{X'}\) is ample.  Thus \(X'\) is klt Fano.  Likewise, the two active
intersection numbers become
\begin{equation}\label{frac:eq:two-line-intersections}
 A_{X'}^rC_{X'}=1,
 \qquad
 B_{X'}A_{X'}^{r-1}C_{X'}=d.
\end{equation}

It remains to prove that \(A_{X'}\) is globally generated.  Unlike in
the one-active case, \(A_{X'}\) is not known to be \(f'\)-ample.
We therefore contract its relative-null directions.  Let
\(h:X'\to W\) be the contraction defined by the semiample divisor
\(A_{X'}+C_{X'}\).  If a curve \(R\) is contracted by \(h\), then
\((A_{X'}+C_{X'})R=0\).  Since both summands are nef, this implies
\(A_{X'}R=C_{X'}R=0\).  Lemma~\ref{lem:descent} therefore gives Cartier
divisors \(A_W,C_W\) on \(W\), with
\[
 A_{X'}=h^*A_W,
 \qquad
 C_{X'}=h^*C_W,
\]
and \(A_W+C_W\) ample.  By
\eqref{frac:eq:two-line-intersections},
\[
 (A_{X'}+C_{X'})^{r+1}
 \geq
 (r+1)A_{X'}^rC_{X'}
 =r+1>0.
\]
Thus \(A_{X'}+C_{X'}\) is big.  Since it is semiample,
\(\dim W=\nu(A_{X'}+C_{X'})=r+1=\dim X'\); hence the contraction \(h\)
is birational.  Since \(C_{X'}=f'^*\OO_{\PP^1}(1)\) is trivial on the
fibres of \(h\), the morphism \(f'\) factors as
\[
 X'\xrightarrow{h}W\xrightarrow{g}\PP^1,
 \qquad
 C_W=g^*\OO_{\PP^1}(1).
\]
Because \(C_W\) is pulled back from the base and \(A_W+C_W\) is ample,
\(A_W\) is \(g\)-ample.

We now absorb the \(B_{X'}\)-term into a boundary.  Since \(B_{X'}\)
is nef Cartier on the klt Fano variety \(X'\), it is semiample.  Choose a
sufficiently divisible \(N\) such that \(|NB_{X'}|\) is base-point-free,
take a general divisor \(D\in|NB_{X'}|\), and put
\(\Delta:=\frac bN D\).  Then \((X',\Delta)\) is klt and
\(\Delta\sim_{\QQ}bB_{X'}\).  Since numerical and \(\QQ\)-linear
equivalence agree on the klt Fano variety \(X'\),
equation~\eqref{frac:eq:two-line} gives
\[
 K_{X'}+\Delta
 \sim_{\QQ}
 h^*\bigl(-(r-t)A_W-(2-\eta)C_W\bigr).
\]
The canonical bundle formula therefore gives an effective
\(\QQ\)-divisor \(\Gamma\) on \(W\) such that \((W,\Gamma)\) is klt and
\begin{equation}\label{frac:eq:two-W-logFano}
 -(K_W+\Gamma)
 \sim_{\QQ}
 (r-t)A_W+(2-\eta)C_W.
\end{equation}
Since \(W\) is klt, the dominant morphism \(g:W\to\PP^1\) is flat, and
its fibres have dimension \(r\).  As \(t<1\), we have
\(r<(r-t)+1\).  Fujino's relative spannedness theorem
\cite[Theorem~1.1]{FujinoSpanned}, applied to
\eqref{frac:eq:two-W-logFano}, shows that \(A_W\) is generated by
relative global sections.  Hence \(A_{X'}=h^*A_W\) is generated
relative to \(f'\).

Now the same relative-to-global argument as in the proof of
Proposition~\ref{frac:prop:one-active} applies: both
\[
 A_{X'}-K_{X'}
 \equiv
 (r+1-t)A_{X'}+bB_{X'}+(2-\eta)C_{X'}
\]
and
\[
 A_{X'}-C_{X'}-K_{X'}
 \equiv
 (r+1-t)A_{X'}+bB_{X'}+(1-\eta)C_{X'}
\]
are ample by \eqref{frac:eq:eta-two-active}.  Hence \(A_{X'}\) is
globally generated.

As in Proposition~\ref{frac:prop:one-active}, cut by \(r-1\) general
members of \(|A_{X'}|\), with no cut when \(r=1\).  Repeated Bertini
and adjunction show that the resulting surface \(S\) is klt.  Writing
\(A_S:=A_{X'}|_S\), \(B_S:=B_{X'}|_S\), and \(C_S:=C_{X'}|_S\),
adjunction gives
\begin{equation}\label{frac:eq:two-surface}
 -K_S\equiv(1-t)A_S+bB_S+(2-\eta)C_S,
 \qquad
 C_S^2=0,
 \qquad
 A_SC_S=1,
 \qquad
 B_SC_S=d.
\end{equation}
All three coefficients are positive, and
\(A_S+B_S+C_S\) is ample; hence \(S\) is Fano.
In particular \(C_S\not\equiv0\).  Since
\(0<\eta<\min\{1,1-t,b\}\) by
\eqref{frac:eq:eta-two-active} and \(b<1\),
Lemma~\ref{frac:lem:surface}~(ii), applied with
\((A_1,A_2)=(A_S,B_S)\) and \((a_1,a_2)=(1-t,b)\), contradicts
\eqref{frac:eq:two-surface}.  Therefore
\begin{equation}\label{frac:eq:inactive-defects-zero-two}
 \eps_j=0
 \qquad(j\in J).
\end{equation}

\smallskip\noindent
\emph{Step 3: reduction to a Picard-number-two core.}
Lemma~\ref{frac:lem:inactive-split} now splits off all inactive factors:
\[
 X\simeq Z\times\prod_{j\in J}\PP^{m_j}.
\]
Since \(s=\rho(X)\) and the only active indices are \(p,q\), we have
\(\rho(Z)=\rho(X)-|J|=2\).  Let \(A\) and \(B\) denote the descended
classes of \(L_q\) and \(L_p\), respectively.  Moreover,
\[
 \dim Z
 =n-\sum_{j\in J}m_j
 =m_p+m_q+1
 =r,
\]
because \(m_p=0\) and \(r=m_q+1\).

The inactive projective-space factors contribute only their top
hyperplane degrees, all equal to one.  Hence
\[
 d_q=L^{m+\ee_q}=A^r=1,
 \qquad
 d_p=L^{m+\ee_p}=BA^{r-1}=d.
\]
Together with the anticanonical decomposition and
\eqref{frac:eq:two-numbers}, this gives
\begin{equation}\label{frac:eq:two-core}
 -K_Z\equiv(r-t)A+bB,
 \qquad
 A^r=1,
 \qquad
 BA^{r-1}=d,
 \qquad
 db=1+t,
 \qquad
 0\leq t<b<1.
\end{equation}
The dimension satisfies $r\geq2$, because a normal projective curve has
Picard number one.

We now identify the contraction defined by \(A\).  Since \(A^r=1>0\)
and \(\dim Z=r\), the nef divisor \(A\) is big.  As \(A\) is semiample,
its contraction \(g:Z\to T\) has
\(\dim T=\nu(A)=r\), so \(g\) is generically finite.  Since \(g\) has
connected fibres, it is birational.  By Lemma~\ref{lem:descent},
\(A\sim g^*H\) for an ample Cartier divisor \(H\) on \(T\), and the
projection formula gives \(H^r=A^r=1\).

For the decomposition \eqref{frac:eq:two-core}, the rounding integers
of the coefficients \(r-t\) and \(b\) are \(r-1\) and \(0\),
respectively.  Thus the vanishing box contains \((z,0)\) for every
\(1\leq z\leq r-1\), and Lemma~\ref{lem:Bm-vanishing} gives
\[
 \chi(Z,\OO_Z(-zA))=0\qquad(1\leq z\leq r-1).
\]
Lemma~\ref{lem:image-P} (ii) therefore identifies
$(T,H)$ with $(\PP^r,\OO(1))$. Thus
\[
 g:Z\longrightarrow\PP^r,
 \qquad A=g^*\OO_{\PP^r}(1),
\]
is birational.

Apply Lemma~\ref{frac:lem:birational} with
\(n=r\), \(\eps=t\), and \(\delta=d\).
All its hypotheses are exactly the numerical relations in
\eqref{frac:eq:two-core}.  If $g$ is an isomorphism, then
$\rho(Z)=1$, contradicting $\rho(Z)=2$.  In the remaining case the
lemma gives
\[
 B\sim d(A-E),
 \qquad
 b=\frac1d,
 \qquad
 t=0,
\]
and \(A-E\) is a nonzero nef Cartier divisor.  Since
\(B\sim d(A-E)\), the summand \(bB\) may be replaced by
\(db(A-E)\).  Here \(t=0\), so \(db=1\), and the coefficient sum
increases from \(r+b\) to \(r+1\), by
\((d-1)b>0\).  This contradicts the optimality of the descended
decomposition on \(Z\).
The assumption of two active indices is therefore impossible.
\end{proof}

\subsection{The rank-one factors}

\begin{proposition}\label{frac:prop:rank-one}
Let $Z$ be a complex klt Fano variety of dimension $r\geq2$ and Picard
number one, let $H$ be a primitive ample Cartier divisor, and let $d\geq3$.
Assume
\[
 -K_Z\equiv\left(r-1+\frac2d\right)H,
 \qquad H^r=d.
\]
Then either
\[
 (Z,H)\simeq\bigl(\PP(1,1,d^{r-1}),\OO(d)\bigr),
\]
or $d=4$ and
\[
 (Z,H)\simeq\bigl(\PP(1,1,1,2^{r-2}),\OO(2)\bigr).
\]
\end{proposition}
\begin{proof}
The assumptions give \(K_Z+(r-1)H\equiv-\frac2dH\).
Since $H$ is ample and $d\geq3$, this divisor is antiample and in
particular is not pseudo-effective.  Thus the adjunction-theoretic
classification in \cite[Theorem~1.4]{Liu} applies to the normal
polarized variety $(Z,H)$.

We explain how the alternatives in that theorem specialize here.  The
cases of projective space and of a quadric have polarized degree one and
two, respectively, whereas
\(H^r=d\geq3.\)
A projective-space bundle over a curve has Picard number at least two,
contrary to $\rho(Z)=1$.  The exceptional polarized surface
$(\PP^2,\OO_{\PP^2}(2))$ has a nonprimitive polarization, whereas $H$
is primitive.  The only remaining alternatives are generalized cones over a
polarized rational normal curve $(\PP^1,\OO_{\PP^1}(a))$, with
$a\geq3$, or generalized cones over $(\PP^2,\OO_{\PP^2}(2))$.
In the first case the polarized degree of the generalized cone is $a$;
since $H^r=d$, we have $a=d$.  In the second case the polarized degree
is four, so necessarily $d=4$.

We use the standard identification of projective cones over Veronese
embeddings with weighted projective spaces; see
\cite[\S1.2.3]{DolgachevWPS} and \cite[p.~304]{BadescuProjective}.
Iterating the same Veronese-ring construction, the $c$-fold generalized
cone over \((\PP^m,\OO_{\PP^m}(q))\) is
\(\bigl(\PP(1^{m+1},q^c),\OO(q)\bigr)\).
Applying this with \((m,q,c)=(1,d,r-1)\) gives
\((Z,H)\simeq\bigl(\PP(1,1,d^{r-1}),\OO(d)\bigr)\).
For the second possibility, applying the same identification with
\((m,q,c)=(2,2,r-2)\) gives
\((Z,H)\simeq\bigl(\PP(1,1,1,2^{r-2}),\OO(2)\bigr)\).
The case $r=2$ is the excluded pair $(\PP^2,\OO(2))$, so in this
alternative $r\geq3$.  These are exactly the two possibilities in the
statement.
\end{proof}

\section{Verification of the examples}\label{sec:converse}

The classification arguments of the preceding sections give the necessary
forms in Theorem~\ref{thm:main}.  We record the converse calculations.

\begin{proposition}\label{prop:examples}
The non-projective-space factors listed in Theorem~\ref{thm:main} are klt
Fano varieties and have the nef complexities asserted there.
\end{proposition}
\begin{proof}
We treat the four types in turn.

Let $Q^d\subset\PP^{d+1}$ be a normal quadric of Picard number one and
let $H=\OO_Q(1)$.  Its quadratic hypersurface singularities are
canonical, hence klt, and adjunction gives
\[
 -K_Q=dH.
\]
Since $H$ is the primitive ample Cartier generator, $\tau_Q=d$ and
$c_Q=d+1-d=1$.  The smooth quadric surface is not among these
Picard-number-one factors: it is $\PP^1\times\PP^1$ and has complexity
zero.

Next let $Y=\Bl_{\PP^\ell}\PP^d$, where $0\leq\ell\leq d-2$.  If
$H$ is the pullback of the hyperplane class, $E$ is the exceptional
divisor, and $M=H-E$, then $H$ and $M$ are the primitive generators of
the nef Cartier monoid.  The blow-up formula gives
\[
 -K_Y=(\ell+2)H+(d-\ell-1)M.
\]
Both coefficients are positive, so $Y$ is Fano, and the displayed
decomposition is optimal.  Hence $\tau_Y=d+1$ and $c_Y=1$.

Let now $Y\subset\PP^a\times\PP^b$, $a,b\geq2$, be an integral
hypersurface of bidegree $(1,1)$.  By
Lemma~\ref{lem:bilinear-singularities}, $Y$ is normal, terminal, and
Gorenstein.  If $H_1,H_2$ are the restrictions of the two hyperplane
classes, adjunction gives
\[
 -K_Y=aH_1+bH_2.
\]
Grothendieck--Lefschetz identifies
$\Pic(Y)=\ZZ H_1\oplus\ZZ H_2$, and the two projections show that
\[
 \Nef(Y)\cap\Lambda_Y
 =\ZZ_{\geq0}H_1+\ZZ_{\geq0}H_2.
\]
Thus $Y$ is Fano, the displayed anticanonical decomposition is optimal,
and $\tau_Y=a+b$, so $c_Y=1$.

Finally consider
\[
 W=\PP(1,1,q^{d-1}),\qquad q\geq2.
\]
The weights are well formed and $W$ has quotient singularities, hence is
$\QQ$-factorial and klt.  Its primitive ample Cartier generator is
$H=\OO_W(q)$.  The dualizing-sheaf formula for weighted projective
spaces \cite[Theorem~3.3.4]{DolgachevWPS} gives
\[
 \omega_W\simeq\OO_W\bigl(-2-(d-1)q\bigr),
\]
and therefore
\[
 -K_W\equiv\left(d-1+\frac2q\right)H.
\]
It follows that $W$ is Fano and, since $H$ is the primitive nef Cartier
generator,
\[
 \tau_W=d-1+\frac2q,
 \qquad
 c_W=2-\frac2q.
\]
For $W=\PP(1,1,1,2^{d-2})$, $d\geq3$, the primitive Cartier generator
is $H=\OO_W(2)$, and the same dualizing-sheaf formula gives
\[
 -K_W\equiv\left(d-\frac12\right)H.
\]
Thus $W$ is klt Fano,
\[
 \tau_W=d-\frac12,
 \qquad
 c_W=\frac32.
\]
\end{proof}

\begin{proof}[Proof of Theorem~\ref{thm:main}]
Statements (i) and (ii) follow from Proposition~\ref{prop:defect} and
Theorem~\ref{thm:gap}.  If $c_X=1$, the splitting
\eqref{eq:splitting} reduces the classification to Picard number at most
two, and Propositions~\ref{prop:rho-one} and \ref{prop:rho-two} give
exactly the three factors in (iii).

If $1<c_X<2$, Corollary~\ref{frac:cor:reduction} gives the fractional
normal form.  Proposition~\ref{frac:prop:two-active} excludes the
two-active case, while Propositions~\ref{frac:prop:one-active} and
\ref{frac:prop:rank-one} split off the inactive projective-space factors
and identify the remaining Picard-number-one factor with one of the two
weighted-projective families in (iv).  Proposition~\ref{prop:examples}
verifies the stated complexities of all non-projective-space factors,
and Lemma~\ref{lem:product-index} shows that adjoining projective-space
factors does not change nef complexity.  This proves the converse
statements as well and completes the proof.
\end{proof}

\begin{proof}[Proof of Corollary~\ref{cor:spectrum}]
Theorem~\ref{thm:main} gives all values below two.  The weighted family
$\PP(1,1,q^{d-1})$ has value $2-2/q$; the value $1$ is obtained for
$q=2$, and the second weighted family has value $3/2$, already obtained
when $q=4$.  Together with the value zero this gives the asserted set.
\end{proof}

\section{The spectrum of nef complexity}\label{sec:spectrum}

We retain $\tau_X$ and $c_X=\dim X+\rho(X)-\tau_X$ from
Definition~\ref{def:nef-complexity} and
$\Lambda_X=\Pic(X)/\!\equiv_{\mathrm{num}}$ from
Section~\ref{sec:preliminaries}.  We construct all rational values at
least two explicitly.

\begin{proposition}\label{prop:spectrum-cones}
Let $q,d,m$ be integers such that
\begin{equation}\label{eq:spectrum-parameters}
 q\geq2,\qquad d\geq1,\qquad m\geq1,\qquad q(m+1)>2d.
\end{equation}
There is a $\QQ$-factorial klt Fano variety $X_{q,d,m}$ of dimension
$m+2$ and Picard number two, with nef Cartier divisors $H,G$ satisfying
\begin{align}
 \Pic(X_{q,d,m})&=\ZZ H\oplus\ZZ G,\label{eq:spectrum-pic}\\
 \Nef(X_{q,d,m})\cap\Lambda_{X_{q,d,m}}
 &=\ZZ_{\geq0}[H]\oplus\ZZ_{\geq0}[G],\label{eq:spectrum-nef}\\
 -K_{X_{q,d,m}}
 &\sim_{\QQ}\left(1+\frac2q\right)H
       +\left(m+1-\frac{2d}{q}\right)G.\label{eq:spectrum-K}
\end{align}
Its total index and nef complexity are
\begin{equation}\label{eq:spectrum-invariants}
 \tau_{X_{q,d,m}}=m+2-\frac{2(d-1)}q,
 \qquad
 c_{X_{q,d,m}}=2+\frac{2(d-1)}q.
\end{equation}
\end{proposition}
\begin{proof}
Put $B=\PP^m$, $T=\PP^1\times B$, let $A,C$ be the pullbacks of
$\OO_{\PP^1}(1)$ and $\OO_B(1)$, and set $L=\OO_T(q,d)=qA+dC$.  With Hartshorne's convention
$\PP_T(\mathcal E)=\operatorname{Proj}_T(\operatorname{Sym}^{\bullet}\mathcal E)$, consider
\[
 p:Y=\PP_T(\OO_T\oplus L)\longrightarrow T,
 \qquad \xi=\OO_Y(1).
\]
The quotient onto $\OO_T$ defines a section
$E\sim\xi-qA-dC$.  If $h:Y\to B$ is the natural projection, then
\[
 \mathcal R:=\bigoplus_{k\ge0}h_*\OO_Y(k\xi),
 \qquad
 \mathcal R_k=\bigoplus_{j=0}^k
 H^0(\PP^1,\OO(jq))\otimes\OO_B(jd).
\]
The usual multiplication maps for $\OO_{\PP^1}(q)$ show that
$\mathcal R$ is generated in degree one.  Set
$X=X_{q,d,m}:=\operatorname{Proj}_B\mathcal R$, let
$H=\OO_X(1)$ and let $G$ be the pullback of $\OO_B(1)$.  The induced
birational morphism $\varphi:Y\to X$ satisfies
$\varphi^*H=\xi$ and $\varphi^*G=C$.

After trivializing $\OO_B(d)$ on an affine open $U\subset B$, one has
$X|_U\simeq U\times\PP(1,1,q)$. Under this identification the
relatively ample line bundle $H$ restricts to
$\OO_{\PP(1,1,q)}(q)$, the primitive Cartier generator.  Hence $X$ is
normal, $\QQ$-factorial and klt, and $\varphi$ contracts exactly $E$
along the $\PP^1$-direction.  Since
$\Pic(Y)=\ZZ\xi\oplus\ZZ A\oplus\ZZ C$ and the contracted curves have
$A$-degree one and zero $\xi$- and $C$-degrees, pullback identifies
\[
 \Pic(X)=\ZZ H\oplus\ZZ G.
\]
Moreover $H$ and $G$ are nef.  Testing a class $aH+cG$ on a fibre
$\PP(1,1,q)$ of $X\to B$ and on a line in the contracted section
$B\subset X$ gives $a,c\ge0$.  Thus $\rho(X)=2$ and
\[
 \Nef(X)=\RR_{\ge0}[H]+\RR_{\ge0}[G],
 \qquad
 \Nef(X)\cap\Lambda_X
 =\ZZ_{\ge0}[H]\oplus\ZZ_{\ge0}[G].
\]

The projective bundle formula gives
\[
 K_Y\sim-2\xi+(q-2)A+(d-m-1)C.
\]
Writing $K_Y=\varphi^*K_X+a_EE$ and comparing the $A$-coefficients gives
$a_E=2/q-1>-1$.  If
$-K_X\sim_{\QQ}xH+yG$, comparison of the $\xi$- and $C$-coefficients
gives
\[
 x=1+\frac2q,
 \qquad
 y=m+1-\frac{2d}{q}.
\]
This proves \eqref{eq:spectrum-K}.  Both coefficients are positive by
\eqref{eq:spectrum-parameters}; since the nef cone is the first quadrant,
Kleiman's criterion \cite[Theorem~1.4.23 (ii)]{Lazarsfeld} shows that
$-K_X$ is ample.

Finally, the nef Cartier monoid is freely generated by $H$ and $G$, so
\eqref{eq:spectrum-K} is optimal and
\[
 \tau_X=m+2-\frac{2(d-1)}q,
 \qquad
 c_X=2+\frac{2(d-1)}q.
\]
This proves \eqref{eq:spectrum-invariants}.
\end{proof}

\begin{corollary}\label{cor:spectrum-realization}
Let $r\in\QQ$ with $r\geq2$. Then $r$ is the nef complexity of a
$\QQ$-factorial klt Fano variety of Picard number two.
\end{corollary}
\begin{proof}
Write $r=2+a/b$ with $a\ge0$ and $b\ge1$, and set
\begin{equation}\label{eq:spectrum-realization-parameters}
 q=2b,
 \qquad d=a+1,
 \qquad m=a+1.
\end{equation}
Then \eqref{eq:spectrum-parameters} holds, and
Proposition~\ref{prop:spectrum-cones} gives a $\QQ$-factorial klt Fano
variety of Picard number two with
\[
 c_X=2+\frac{2(d-1)}q=2+\frac ab=r.
\]
\end{proof}

\begin{proof}[Proof of Theorem~\ref{thm:full-spectrum}]
Proposition~\ref{prop:LP} shows that the total index of every klt Fano
variety is rational; hence every nef complexity is rational.  Theorem~\ref{thm:main}
and Corollary~\ref{cor:spectrum} determine all values strictly below two.
Conversely, Corollary~\ref{cor:spectrum-realization} constructs a
$\QQ$-factorial klt Fano variety of Picard number two for every rational
number $r\geq2$.  Combining these two statements gives the asserted full
spectrum.
\end{proof}

\begin{remark}
The dimension in \eqref{eq:spectrum-realization-parameters} is $a+3$ and
is not fixed as the value varies. All varieties in this construction are
singular because $q\geq2$.
\end{remark}

\section*{Acknowledgments}
The author thanks Professor Yoshinori Gongyo for listening to a presentation of
an alternative proof of $c_X\geq0$, and Professor Joaqu\'in Moraga for discussions
on related questions. This work was partially supported by JSPS KAKENHI
Grant Number 25K06940.

\end{document}